\documentclass[a4paper]{article}
\usepackage{amsmath, amssymb}
\usepackage{stmaryrd}
\usepackage{amsthm,cite}
\usepackage{amsfonts}
\usepackage{graphicx,subfigure}
\usepackage{epsfig,color,pifont,enumerate}
\usepackage{setspace}
\graphicspath{{img/}}

\usepackage{geometry}
\newcommand{\real}{\mathbb{R}}

\renewcommand{\phi}{\varphi}
\renewcommand{\epsilon}{\varepsilon}
\renewcommand{\kappa}{\varkappa}
\newcommand{\sgn}{\text{\bf sgn}\,}
\newcommand{\trs}{\top}
\newcommand{\br}{\mathbb{R}}
 \newcommand{\bldr}{\boldsymbol{r}}

\newtheorem{theorem}{Theorem}[section]

\newtheorem{remark}{Remark}[section]
\newtheorem{assumption}{Assumption}[section]
\newtheorem{lemma}{Lemma}[section]
\newtheorem{cor}{Corollary}[section]
\newtheorem{proposition}{Proposition}[section]

\newcommand{\so}{\text{\scriptsize $\mathcal{O}$}}
\newcommand{\spr}[2]{\left\langle #1; #2 \right\rangle}

\newcommand{\ov}[1]{\overline{#1}}

\newcommand{\ve}{\varepsilon}

\newcommand{\dd}[1]{#1^{\prime\prime}}

\newcommand{\pf}{{\bf Proof:}\;}
\newcommand{\epf}{$\qquad \Box$}
\newcommand{\dist}[2]{\text{\bf dist}\, \left[ #1; #2 \right]}
\newcommand{\pp}{{\prime\prime}}
\title{Proofs of the Technical Results Justifying an Algorithm of Extremum Seeking Navigation in Dynamic Environmental Fields}
\date{}
\author{Alexey S. Matveev, Michael C. Hoy, and Andrey V.Savkin}
\begin{document}
\maketitle
\section{Introduction}
The problem of extremum seeking navigation is concerned with driving a mobile
robot to the maximizer of an environmental field. The field profile is unknown a priori and should be explored on-line based on point-wise measurements of the field. This scenario arises in a variety of missions, such as environmental studies, tracking a moving target based on a single signal decaying away from the target, seeking sources of radiation emission or pollutant leakage, etc.
The discussed problem falls into the general framework of extremum seeking control, which should ensure that a dynamic plant is directed to and subsequently tracks an optimal operating point based on on-line evaluation of the current performance, whose profile is unknown a priori.
\par
For recent surveys on extremum seeking control methods and algorithms, we refer the reader to \cite{KoRu08,TaMoMaNeMa10,DoPeGu11}.
A good deal of related research was concerned with gradient climbing based on direct on-line gradient estimation \cite{BaChSp00,TePo01,LiDu04,AtNyMiPa12,ZhOr12}. This approach is especially beneficial for mobile sensor networks thanks to collaborative field
measurements in many locations and data exchange \cite{PoNe96,OgFiLe04,BiAr07,Cor07,BaLe02,GaPa04,MoCa10,PePa11}.
However even in this scenario, data exchange degradation due to e.g., communication constraints may require each robot in the team to operate autonomously over  considerable time. Similar algorithms can be basically used for a single robot equipped with several sensors that are distant enough from each other and thus provide the field values at several essentially diverse locations. In any case,  multiple vehicle/sensor scenario means
complicated and costly hardware. This paper is focused on another case, where at any time, the robot has access to only a single field measurement at its current location.
\par
The lack of multiple sensor data can be compensated via exploring multiple nearby locations by "dithering" the position of the single sensor during special maneuvers, which may be excited by probing high-frequency sinusoidal \cite{BuYoBrSi96,ZhArGhSiKr07,CoKr09,CoSiGhKr09} or stochastic \cite{LiKr10} inputs. An approach that is similar in spirit is extremum seeking by means of many robots performing biased random walks \cite{MeHeAs08} or by two robots with access to relative positions of each other and rotational actuation \cite{ElBr12}.
These methods rely, either implicitly or explicitly, on systematic sideways exploration maneuvers to collect rich enough data. Another approach limits the field gradient information to only the time-derivative of the measured field value obtained via e.g., numerical differentiation, \cite{BaLe02,MaSaTe08,BarBail08,MaTeSa11} and partly employs switching controllers \cite{MaSaTe08,MaTeSa11}.
These give rise to concerns about amplification of the measurement noises and chattering, respectively.
Though the respective potential detrimental consequences can be successfully avoided in the scenarios
examined in these papers, anyhow the need to ensure this puts extra burden on controller parameters tuning.
Adaptive extremum seeking approach \cite{GaZh03} assumes that the field is known up to finitely many uncertain steady parameters, which may exceed the real level of knowledge in some applications.
\par
The previous extremum seeking research dealt mainly with steady fields. However in the real world, environmental fields are almost never steady and often cannot be well approximated by steady fields. The problem of extremum seeking in dynamic fields is also concerned with navigation and guidance of a mobile robot towards an unknowingly maneuvering target based on a single measurement that decays as the sensor goes away from the target, like the strength of the infrared, acoustic, or electromagnetic signal, or minus the distance to the target. Such navigation is of interest in many areas \cite{ADB04,GS04,MTS11ronly}; it carries a potential to reduce the hardware complexity and cost and improve target pursuit reliability. A solution for such problem in the very special case of the unsteady field --- minus the distance to an unknownly moving Dubins-like target --- was proposed and justified in \cite{MTS11ronly}. However the results of \cite{MTS11ronly} are not applicable to more general dynamic fields.
\par
Unlike the previous research, this paper deals with generic dynamic fields. In this context, it justifies a new kinematic control paradigm that offers to keep the velocity orientation angle proportional to the discrepancy between the field value and a given linear ascending function of time, as opposed to conventionally trying to align the velocity vector with the gradient.
This control law is free from evaluation of any field-derivative data, uses only finite gains instead of switching control, and demands only minor memory and computational robot's capacities, being reactive in its nature.
We provide a mathematically rigorous evidence of convergence of this
control law in the case of a generic dynamic field. In doing so, its is shown that the closed-loop system is prone to monotonic, non-oscillatory behavior during the transient to the field maximizer provided that the controller parameters are properly tuned.
We offer recommendations on the choice of these parameters under which the robot inevitably reaches the desired vicinity of the moving field maximizer in a finite time and remains there afterwards.
\par
 An algorithm similar to ours can be found in \cite{BarBail08}, which however uses the estimated time-derivative of the measurement. Another difference is that in \cite{BarBail08}, only a steady harmonic field was examined, the performance during the transient and the behavior after reaching a vicinity of the maximizer were not addressed even for general harmonic fields, and the convergence conditions were partly implicit by giving no explicit bound on some entities that were assumed sufficiently large. The focus of this paper is on general dynamic fields, and we offer a study of the entire maneuver with explicit conditions for convergence.
\par
The body of the paper is organized as follows. Section~\ref{sec1}
presents the problem setup and control law. Section~\ref{sec.linear} discloses the closed-loop behavior in a simple but instructive case of a linear field. The main theoretical results are given in Section~\ref{sec.maxx}, where a generic dynamic field is considered.
These results are illustrated by examples in Section~\ref{sec.exampl}, and are confirmed and supplemented via simulation tests and experiments with a real wheeled robot in Sections~\ref{sec.simtest} and \ref{sec.exper}, respectively. Section~\ref{sec.concl} offers brief conclusions.
The proofs of all theoretical results are given in three appendices.
\par
\par
The extended introduction and discussion of the proposed control law are given in the paper submitted by the authors to the IEEE Transactions on Control Systems Technology.
This text basically contains the proofs of the technical facts underlying justification of the convergence at performance of the proposed algorithm in that paper, which were not included into it due to the length limitations. To make the current text logically consistent, we reproduce the problem statement and notations.

\section{Problem Setup and the Control Algorithm}\label{sec1}
We consider a point-wise robot traveling in a Cartesian plane with the absolute
coordinates $x$ and $y$.
The robot is controlled by the time-varying linear velocity $\vec{v}$ whose magnitude does not exceed a given constant $\ov{v}$. The plane
hosts an unknown scalar time-varying field $D(t,\boldsymbol{r}) \in \br$, where $t$ is time and
$\boldsymbol{r}:= (x,y)^\trs$. The objective is to drive
the robot to the point $\bldr^0(t)$ where
$D(t,\bldr)$ attains its maximum over $\bldr$ and then to keep it in a vicinity of
$\bldr^0(t)$, thus displaying the approximate location of
$\bldr^0(t)$. The on-board control system  has access only to the field
value $d(t):= D[t,\bldr(t)]$ at the robot current location $\bldr(t) = [x(t),y(t)]^\trs$.
No data about the derivatives of $D$ are available; in particular, the robot is aware of neither the gradient of $D(\cdot)$ nor the time-derivative $\dot{d}$ of the measurement $d$.
\par
The kinematic model of the robot is as follows:
\begin{equation}
\label{1}
\dot{\bldr} = \vec{v}, \qquad \bldr (0) = \bldr_{\text{in}}\qquad \|\vec{v}\| \leq \ov{v},
\end{equation}
where $\|\cdot \|$ is the Euclidian norm.
The problem
is to design a controller that drives the
robot into the desired vicinity $V_\star(t)$
of the time-varying maximizer $\boldsymbol{r}^0(t)$ in a finite time $t_0$  and then keeps the robot within $V_\star(t)$.
\par
In this paper, we examine the following control algorithm:
\begin{equation}
\label{c.a}
\vec{v}(t) = \ov{v} \vec{e}\Big\{\mu\big[d(t) - d_0(t) \big] \Big\}, \quad
\text{where} \quad d_0(t) := \nu t+d_\ast \quad \text{and} \quad
\vec{e}(\theta) := \left( \begin{smallmatrix} \cos \theta \\ \sin \theta \end{smallmatrix}\right)
\end{equation}
is the unit vector with the heading angle $\theta$ and $\nu>0 , \mu >0, d_\ast \in \br$ are parameters of the controller.
The control law \eqref {c.a} keeps the robot heading proportional to the discrepancy $d(t) - d_0(t)$ between the current field value $d(t)$ and the linear reference signal $d_0(\cdot)$ that ascends at the rate of $\nu$ and is interpreted as a desired or requested behavior of the field over the robot's trajectory.

\section{Sample Behavior in a Linear Field}
\setcounter{equation}{0}
\label{sec.linear}
We start our study of the proposed control law from analysis of its effect in steady linear fields
\begin{equation}
\label{linear.field}
D(x,y) = n \big( x \cos \varphi + y \sin \varphi \big)+D_0, \qquad \text{where} \quad n>0, \varphi, \; \text{and}\;D_0 \in \br
\end{equation}
are given. This is instructive since any smooth dynamic field is well approximated by a linear steady field in a sufficiently small (and sometimes not so small) region of space-time. So such analysis discloses basic behavioral primitives that underlay, more or less, the closed-loop performance in generic fields. We first assume that the requested field ascending rate $\nu = \dot{d}_0$ is feasible. For vehicles traveling at speeds no greater than $\ov{v}$, feasible rates $\dot{d}$ do not exceed $\ov{v}n$.
So we first assume that $\nu < \ov{v}n$; the opposite case will be addressed in Lemma~\ref{lem.typical}.
\par
We first note that perfect tracking $d(t) \equiv d_0(t)$ of the reference signal is not urgent: for gradient climbing, it suffices to ensure ascend of the field value $d$ at the requested rate $\nu$. Since
\begin{equation}
\label{dot.d}
\dot{d} \overset{\text{\eqref{linear.field}}}{=\!=} n \big( \dot{x} \cos \varphi + \dot{y} \sin \varphi \big) \overset{\text{\eqref{1}}}{=\!=}
n v \big( \cos \theta \cos \varphi + \sin\theta \sin \varphi \big) \overset{\text{\eqref{c.a}}}{=\!=} n \ov{v} \cos(\theta-\varphi),
\end{equation}
this is possible only if the robot moves over a straight line subtending an angle of $\theta_\dagger = \varphi \pm \arccos \frac{\nu}{\ov{v}n} + 2 \pi k$ with the $x$-axis, where $k$ is an integer. In turn, such a motion conforms to the control law \eqref{c.a} if and only if
\begin{multline*}
\mu^{-1}\theta_\dagger = d(t) - d_0(t) =
n \left[ \left( x_{\text{in}} + t \ov{v} \cos \theta_\dagger\right) \cos \varphi + \left( y_{\text{in}} + t \ov{v} \sin \theta_\dagger\right) \sin \varphi \right] - \nu t - d_\ast
\\
= \underbrace{\left[ n \ov{v} \cos (\varphi- \theta_\dagger) - \nu \right]}_{=0}\cdot t +
D(\bldr_{\text{in}})  - d_\ast  = D(\bldr_{\text{in}})  - d_\ast
\end{multline*}
or, in other words,
\begin{equation}
\label{lin.equat}
\mu \left[ D(\bldr_{\text{in}})  - d_\ast \right] = \theta_\dagger:= \varphi \pm \theta_\ast + 2 \pi k, \qquad \theta_\ast:= \arccos \frac{\nu}{\ov{v}n}.
\end{equation}
Thus any line $L$ subtending an angle of $\theta_\ast$ with the field gradient accommodates infinitely many ``perfect''
 trajectories of the closed-loop system; their characteristic feature is that the requested field growth rate $\nu$ is tracked perfectly $\dot{d} \equiv \nu$. These processes are enumerated
by $k=0,\pm 1, \ldots$; the $k$th of them is at the speed $\ov{v}$ along the line $L$ and from the unique point $\bldr_{\text{in}} \in L$ that solves the linear equation \eqref{lin.equat}; see Figs.~{\rm \ref{fig,series}(a,b)}. The starting points related to various $k$ are equally spaced and separated by the distance $2 \pi/(\mu n)$.
\begin{figure}[h]
\centering
\subfigure[]{\scalebox{0.3}{\includegraphics{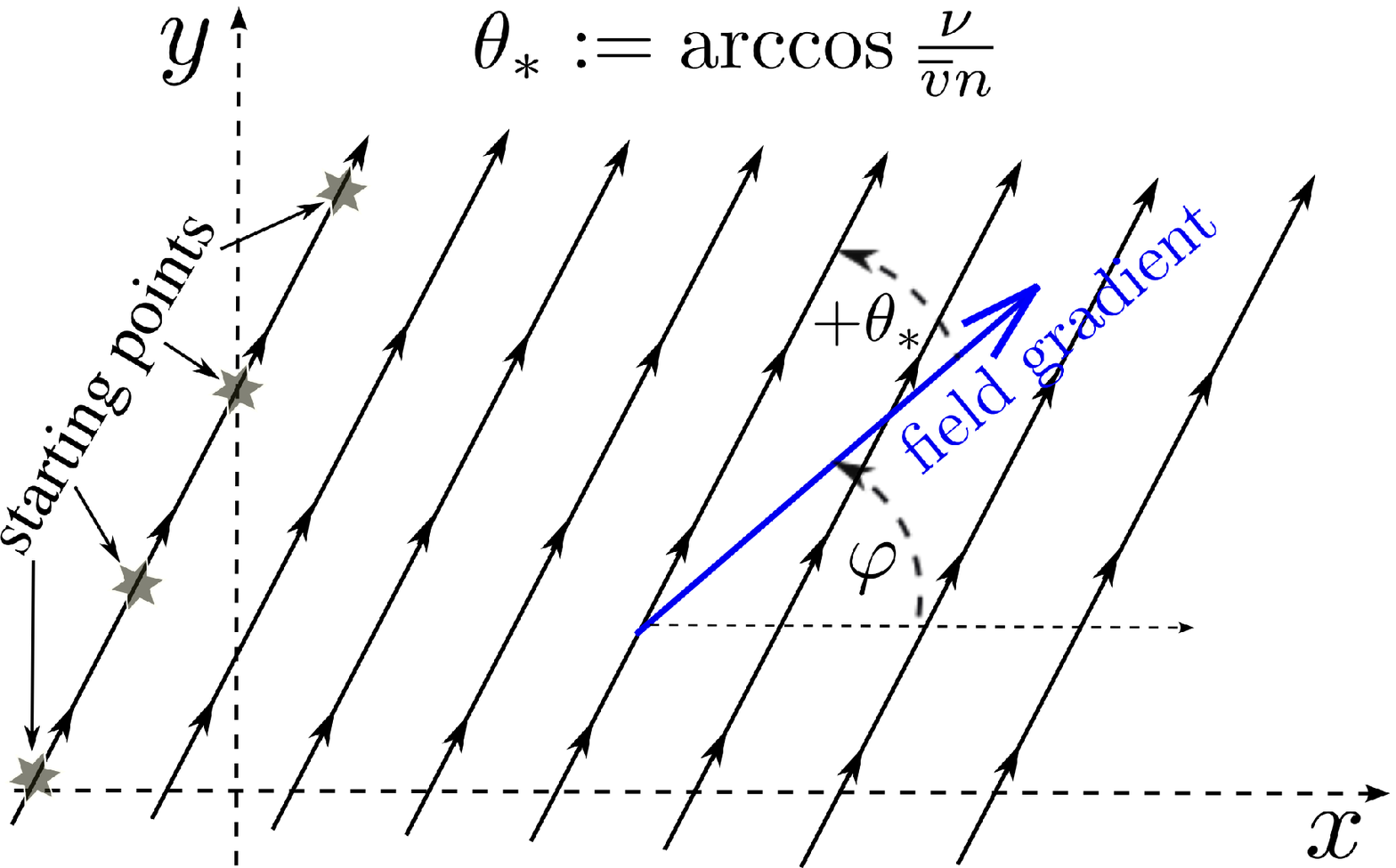}}}
\hfil
\subfigure[]{\scalebox{0.3}{\includegraphics{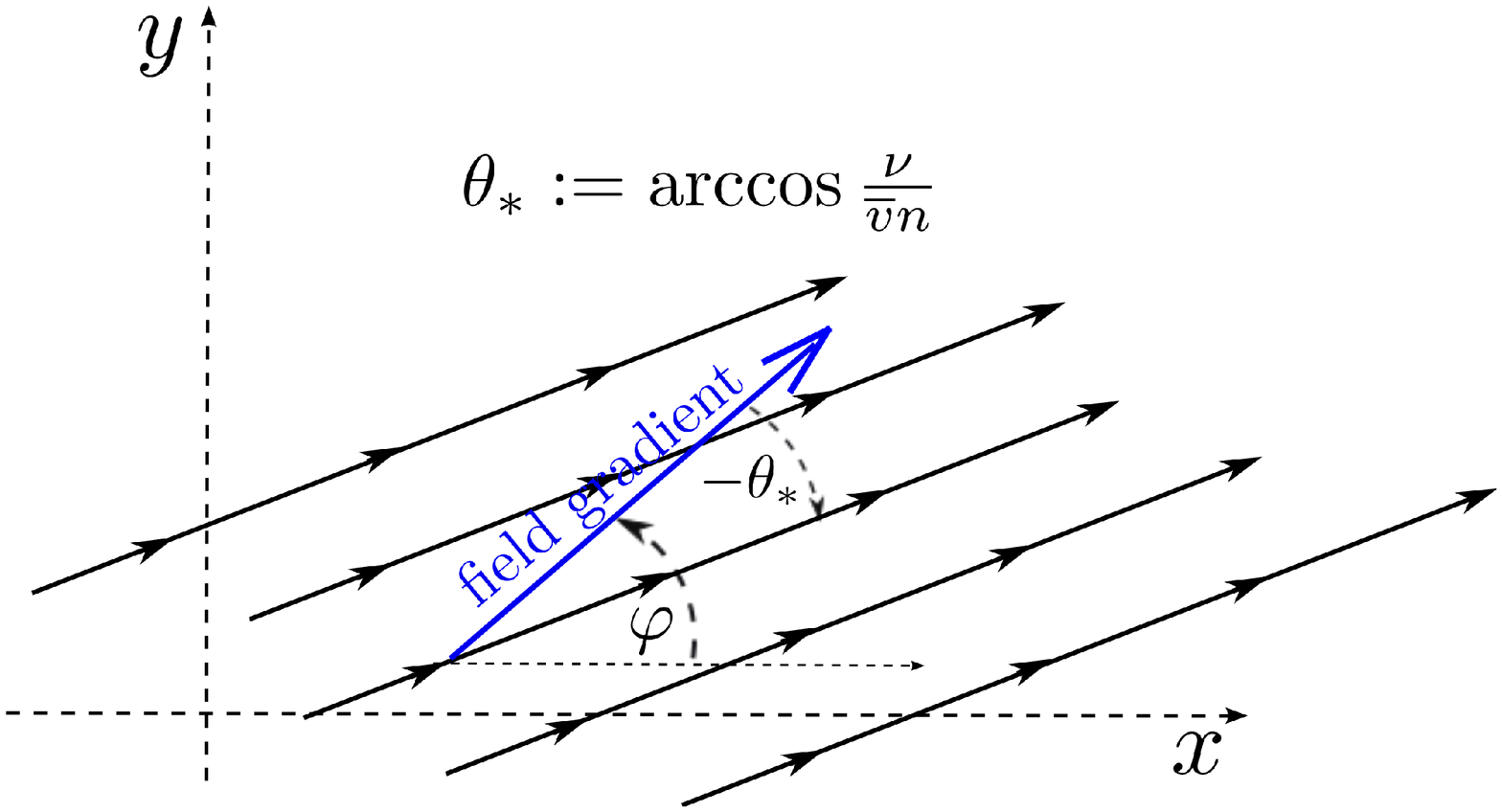}}}
\hfil
\subfigure[]{\scalebox{0.22}{\includegraphics{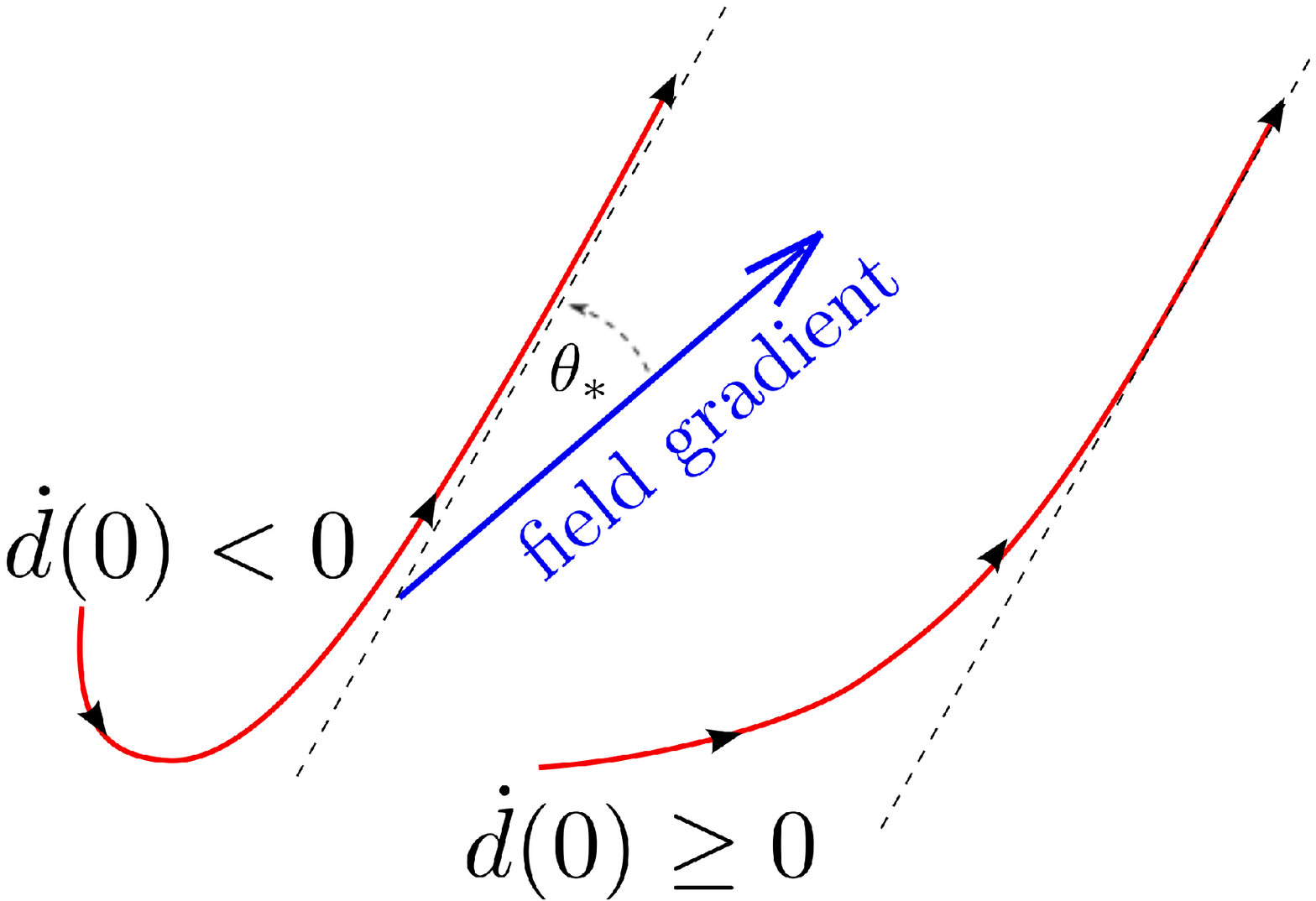}}}
\caption{(a,b) ``Perfect'' closed-loop trajectories (a) $\theta = \varphi + \theta_\ast$ and (b) $\theta = \varphi - \theta_\ast$; (c) Convergence to a ``perfect'' trajectory.}
\label{fig,series}
\end{figure}
\par
\begin{figure}
\centering
\subfigure[]{\scalebox{0.3}{\includegraphics{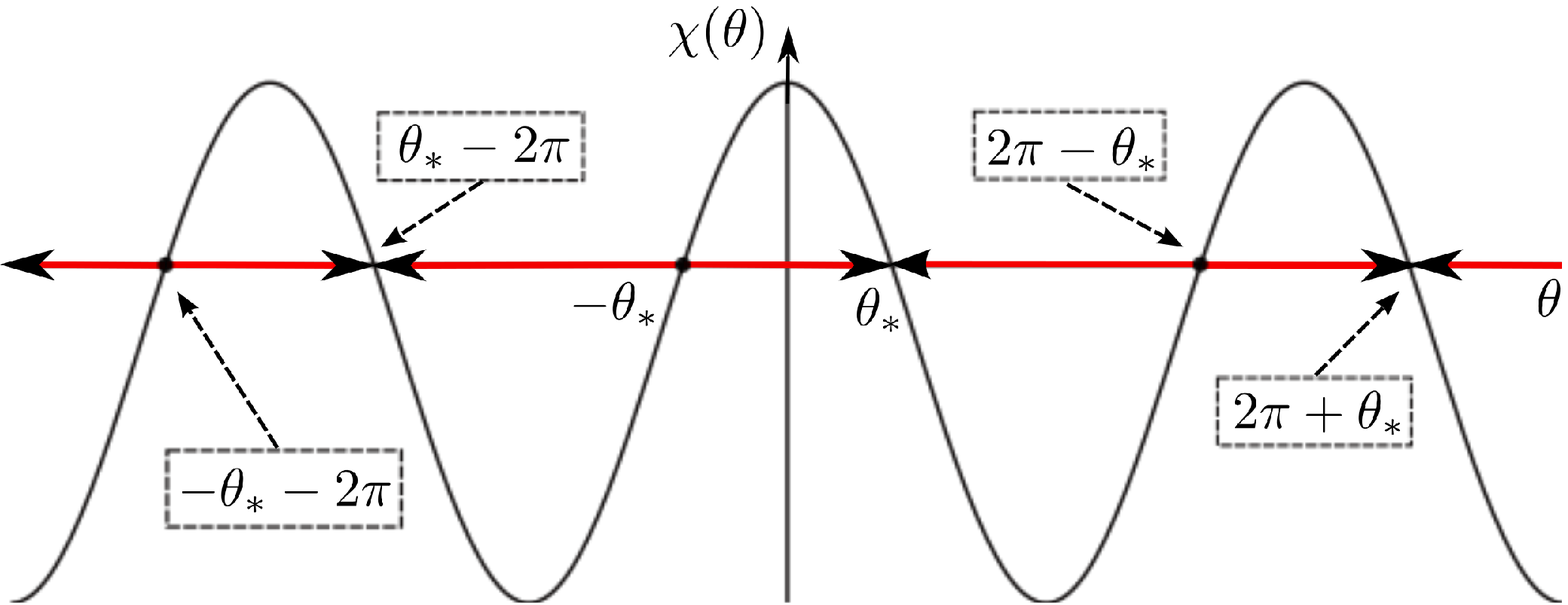}}}
\hfil
\subfigure[]{\scalebox{0.15}{\includegraphics{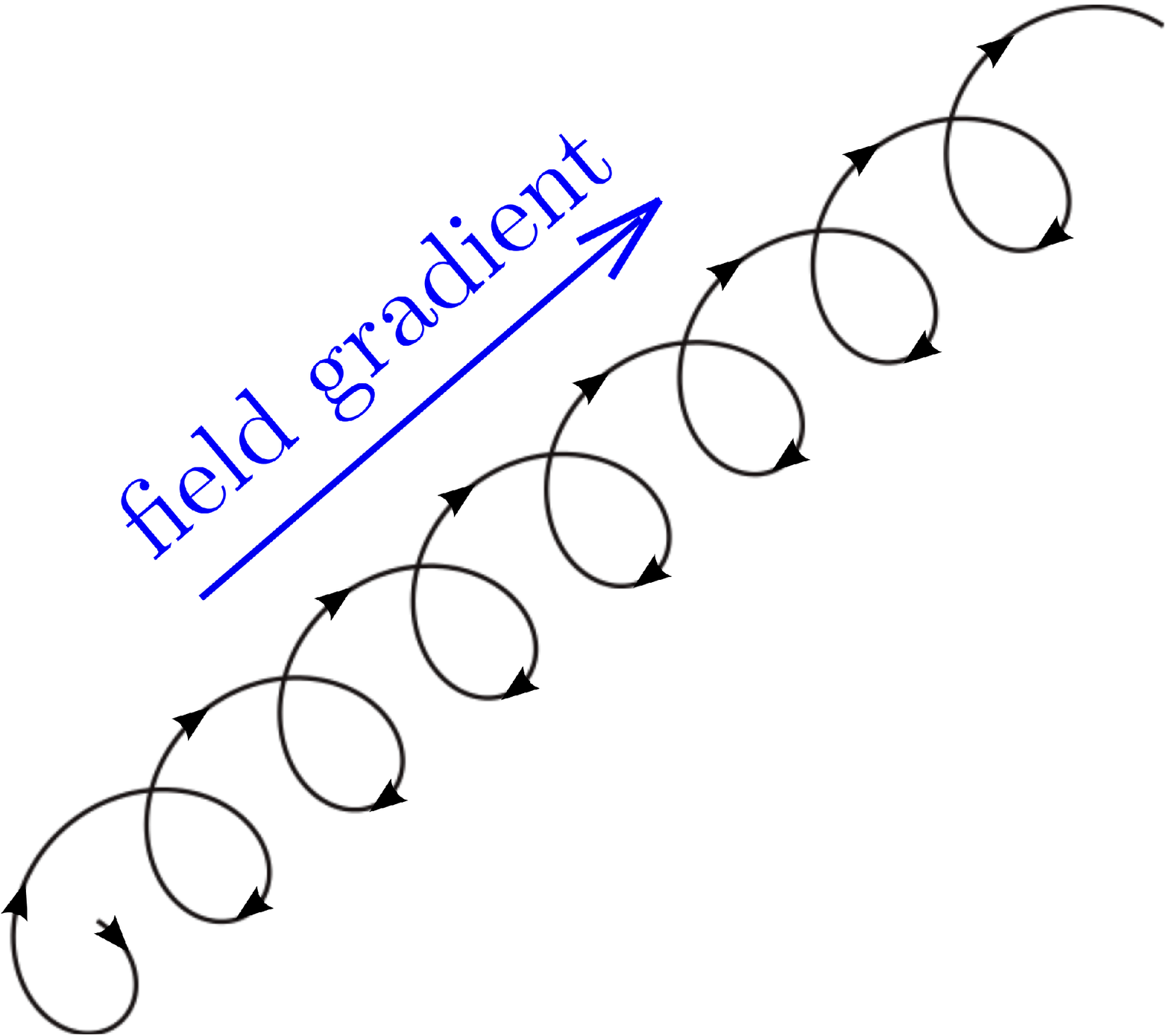}}}
\caption{\small (a) Function $\chi(\cdot)$ and the angular behavior of the robot; (b) Typical closed-loop behavior for $\nu > \ov{v} n$.}
\label{chitheta.fig}
\end{figure}

The following theorem summarizes the above analysis and addresses stability of ``perfect'' processes.
It also shows that the properly tuned controller makes the system pre-disposed to monotonic, unfluctuating behavior.
\begin{theorem}
\label{th.linear}
Let the system \eqref{1} be driven by the controller \eqref{c.a} in the steady linear field \eqref{linear.field}, and the requested field growth rate be feasible $\nu <  \ov{v} n$.
Depending on the initial state $\bldr_{\text{in}}$, all trajectories are divided into two categories:
\begin{enumerate}[{\bf i)}]
\item Let $\bldr_{\text{in}}$ be such that
\eqref{lin.equat} holds.
Then the robot moves at the speed $\ov{v}$ along a straight line subtending the angle $\theta_\ast$ from \eqref{lin.equat} with the field gradient. The field value $d$ constantly grows at the requested rate: $\dot{d} \equiv \nu$;
    \item If \eqref{lin.equat} does not hold,
    the trajectory of the robot, including the velocity orientation,
    exponentially converges to a ``perfect'' trajectory described in {\bf i)} and Fig.~{\rm \ref{fig,series}(a)}.
Meanwhile, the field value $d$ eventually ascends and ultimately approaches the requested ascending rate $\dot{d}(t) \to \nu$ as $t \to \infty$. As $t$ runs from $0$ to $\infty$, the velocity orientation $\theta$ monotonically sweeps through an angle that is less than $2 \pi$, like in Fig.~{\rm \ref{fig,series}(c)}.
\end{enumerate}
\end{theorem}
\pf
With a proper rotation of the coordinate frame in mind, it can be assumed that $\varphi =0$ in \eqref{linear.field}.
\par
i) The claim summarizes the study of ``perfect'' trajectories carried out in the foregoing.
\par
ii) We start with observing that
\begin{equation}
\label{eq.chi}
\dot{\theta} \overset{\text{\eqref{c.a}}}{=\!=} \mu \big[ \dot{d} - \dot{d}_0\big] \overset{\text{\eqref{dot.d}}}{=\!=}
\mu \big[ n \ov{v} \cos \theta - \nu \big] \overset{\text{\eqref{lin.equat}}}{=\!=} \mu \ov{v} n \chi(\theta), \quad \text{where} \quad \chi(\theta) := \cos\theta - \cos \theta_\ast.
\end{equation}
Here $\chi(\cdot)$ has infinitely many roots $\theta_\pm (k) = \pm \theta_\ast + 2 \pi k, k=0,\pm 1, \ldots$ (see Fig.~\ref{chitheta.fig}(a)) and $\pm \chi^\prime [\theta_\pm(k)] <0$.
Hence $\theta_+ (k)$ and $\theta_- (k)$
are locally stable and unstable equilibria, respectively, and any solution $\theta(\cdot) \not\equiv \theta_- (k) \; \forall k $ of \eqref{eq.chi} exponentially and monotonically converges to $\theta_+ (k)$ with some integer $k$, moving within an interval lying within two adjacent equilibria \cite{Hart82}. It remains to note that this convergence $\theta(t) \to \theta_+(k)$ implies the exponential convergence $\bldr(t) - \bldr_+(t) \to 0$ of the trajectories $\bldr(t) = \ov{v} \int_0^t \vec{e}[\theta(s)]\; ds + \bldr_{\text{in}}$ and $\bldr_+(t) = \ov{v} \int_0^t \vec{e}[\theta_+(k)]\; ds + \bldr_{\ast} $, where $\bldr_\ast:= \bldr_{\text{in}} + \ov{v} \int_0^\infty \big\{ \vec{e}[\theta(s)] - \vec{e}[\theta_+(k)]\big\} \; ds $ satisfies \eqref{lin.equat} (with $\bldr_{\text{in}} := \bldr_\ast, \varphi=0, \theta_\dagger := \theta_+(k)$):
\begin{multline*}
\mu \big[ D(\bldr_\ast) - d_\ast\big] = \mu \big[ D(\bldr_{\text{in}}) - d_\ast\big] + \mu \ov{v}n \int_0^\infty \big[\cos \theta(s) - \cos \theta_+(k)\big] \; ds
\\
= \mu \big[ d(0) - d_\ast\big] + \mu \ov{v}n \int_0^\infty \big[\cos \theta(s) - \cos \theta_\ast\big] \; ds
\overset{\text{\eqref{c.a}, \eqref{eq.chi}}}{=\!=\!=\!=\!=\!=} \theta(0) + \int_0^\infty \dot{\theta} (s) \; ds = \lim_{t \to \infty} \theta(t) = \theta_+(k). \qquad \text{\epf}
\end{multline*}
\par
Claim ii) implies that any ``perfect'' trajectory described in {\bf i)} and Fig.~{\rm \ref{fig,series}(b)} is unstable, contrary to those from Fig.~{\rm \ref{fig,series}(a)}.
By the last sentence from ii), the robot does not oscillate and does not make full turns since oscillation violates monotonic evolution of $\theta$, whereas a full turn means that $\theta$ runs the length of $2 \pi$.
This monotonic pattern is antipodal to perpetual oscillations caused by the controllers from \cite{BuYoBrSi96,ZhArGhSiKr07,MeHeAs08,CoKr09,CoSiGhKr09,LiKr10,ElBr12} even in linear fields.
\par
By retracing the arguments from the proof of Theorem~\ref{th.linear}, it is easy to see that
in the marginal case of the maximal feasible field-rate request $\nu = \ov{v} n$, the robot still successfully achieves this rate $\dot{d}(t) \to \nu$ as $t \to \infty$, with retaining the monotonic behavior described by the last sentence from ii).
\par
As for the case where unrealistic field growth rate is requested $\nu > \ov{v}n$, the following lemma shows that even in this case, motion in the gradient direction is the overall average. However the price for unrealistic request is dismissal of the monotonic pattern from Fig.~\ref{fig,series}(c): instead, the robot spirals along a spring-like path, like in Fig.~\ref{chitheta.fig}(b). Such path can be viewed as the result of a periodic motion along a simple closed curve in a clockwise direction, which motion holds in a reference frame that is translating at a constant speed $w >0$ in the gradient direction $\vec{n}_0:= (\cos \varphi, \sin \varphi)^\trs$.
\begin{lemma}
\label{lem.typical}
Suppose that the requested field growth rate is unrealistic $\nu >  \ov{v} n$.
Then the closed-loop trajectory $\bldr(t)$, including the velocity orientation angle $\theta(t)$, has the following form
\begin{equation}
    \label{cyclic.motion}
    \bldr(t) = \bldr_{\circlearrowleft}(t-t_0) + \bldr_\ast +  (t-t_0) w \vec{n}_0  , \quad \theta(t) = \theta_{\circlearrowleft}(t - t_0) - \frac{2 \pi}{\tau}(t-t_0).
\end{equation}
Here $\theta_{\circlearrowleft}(\cdot) : \br \to \br$ and
$\bldr_{\circlearrowleft} (\cdot) : \br \to \br^2$ are periodic functions with a common period $\tau>0$, which along with the speed $w>0$, do not depend on the robot's initial location $\bldr_{\text{in}}$, whereas the constants $t_0 \in \br$ and $\bldr_\ast \in \br^2$ are determined by $\bldr_{\text{in}}$. The equation $\bldr = \bldr_{\circlearrowleft} (t)$ describes a motion along a simple closed curve in the clockwise direction, and the robot's velocity constantly rotates clockwise: $\dot{\theta} <0$.
\end{lemma}
\pf
Like in the proof of Theorem~\ref{th.linear}, we assume, without any loss of generality, that $\varphi=0$ in \eqref{linear.field}.
We also put $\theta_\star:= \arccos \frac{\ov{v}n}{\nu}$. It is easy to see that now \eqref{eq.chi} should be replaced by
\begin{equation}
\label{rhside}
\dot{\theta} = \mu \big[ \ov{v} n \cos \theta  - \nu\big]
= -\mu \nu \big[ 1 - \cos \theta  \cos \theta_\star \big] \leq -\mu \nu \big[ 1 - \cos \theta_\star \big] <0 .
\end{equation}
So any solution $\theta(t)$ goes from $+\infty$ to $-\infty$ as $t$ runs from $- \infty$ to $\infty$ and is a time-shift $\theta(t) = \theta_\lozenge(t - t_0)$ of the unique solution $\theta_\lozenge(\cdot)$ with $\theta_\lozenge(0)=0$. Here $\theta_\lozenge(-t) = - \theta_\lozenge(t)$ since $\theta(t):=-\theta_\lozenge(-t)$ clearly solves \eqref{rhside} and $\theta(0)=0$. Furthermore, $$
\theta(t+\tau) = \theta(t)- 2\pi, \quad \text{where} \quad \tau := \frac{1}{\mu \nu}\int_{-\pi}^\pi \frac{d \theta}{1 - \cos \theta  \cos \theta_\star } >0.
$$
So the function $\theta_{\circlearrowleft}(t):= \theta_\lozenge(t)+\frac{2 \pi}{\tau}t$ is $\tau$-periodic. Hence the second relation from \eqref{cyclic.motion} and the last claim of the lemma do hold.
Taken $\theta$ as an independent variable, we have by \eqref{1} and \eqref{rhside},
\begin{equation*}
\frac{d x}{d \theta} = - \frac{\ov{v}}{\mu \nu} \frac{\cos \theta}{1 - \cos \theta  \cos \theta_\star},
\qquad
\frac{d y}{d \theta} = - \frac{\ov{v}}{\mu \nu} \frac{\sin \theta}{1 - \cos \theta  \cos \theta_\star}.
\end{equation*}
Let $x_\star(\theta)$ and $y_\star(\theta)$ be the respective solutions with $x_\star(0) = y_\star(0)=0$. Since $\sin$ and $\cos$ are odd and even, respectively,
\begin{eqnarray}
\label{inters.that}
\int_{-\pi}^\pi \frac{\sin \theta}{1 - \cos \theta  \cos \theta_\star} \; d\theta =0 \Rightarrow y_\star(\theta+2 \pi) = y_\star(\theta); \quad y_\star(- \theta) = y_\star(\theta), \dot{y}_\star(\theta) <0 \; \forall \theta \in (0,\pi);
\\
\nonumber
\int_{-\pi}^\pi \frac{\cos \theta}{1 - \cos \theta  \cos \theta_\star} \; d\theta =
\int_{-\frac{\pi}{2}}^{\frac{3}{2}\pi} \frac{\cos \theta}{1 - \cos \theta  \cos \theta_\star} \; d\theta =
\int_{-\frac{\pi}{2}}^{\frac{\pi}{2}} \frac{\cos \theta}{1 - \cos \theta  \cos \theta_\star} \; d\theta + \int_{-\frac{\pi}{2}}^{\frac{\pi}{2}} \frac{\cos (\theta+\pi)}{1 - \cos (\theta +\pi) \cos \theta_\star} \; d\theta
\\
\nonumber
=
\int_{-\frac{\pi}{2}}^{\frac{\pi}{2}} \frac{\cos \theta}{1 - \cos \theta  \cos \theta_\star} \; d\theta - \int_{-\frac{\pi}{2}}^{\frac{\pi}{2}} \frac{\cos \theta}{1 +\cos \theta  \cos \theta_\star} \; d\theta
 =
2\int_{-\frac{\pi}{2}}^{\frac{\pi}{2}} \frac{\cos^2 \theta \cos \theta_\star}{(1 - \cos \theta  \cos \theta_\star)(1 +\cos \theta  \cos \theta_\star)} \; d\theta =: L>0
\\
\label{x.even}
\Rightarrow
x_\star(\theta+2 \pi) = x_\star(\theta) - \frac{\ov{v}}{\mu \nu} L ; \qquad x_\star(\theta) = - x_\star(- \theta).
\end{eqnarray}
Hence the function $x_{\star\star}(\theta):= x_\star(\theta)+ \frac{\ov{v}L}{2 \pi \mu \nu} \theta$ is $2 \pi$-periodic.
So $x_{\circlearrowleft}(t):= x_{\star\star}[\theta_\lozenge(t)] - \frac{\ov{v}L}{2 \pi \mu \nu} \theta_{\circlearrowleft}(t) = x_\star[\theta_\lozenge(t)]- \frac{\ov{v}L}{\mu\nu \tau}t$ and $y_{\circlearrowleft}(t) := y_\star[\theta_\lozenge(t)]$
are $\tau$-periodic and the first relation from \eqref{cyclic.motion} holds with $w:= \frac{\ov{v}L}{\mu\nu \tau}$. It remains to show that the curve with the parametric representation $x = x_{\circlearrowleft}(t), y = y_{\circlearrowleft}(t)$ does not intersect itself during its period $t \in (-\tau/2, \tau/2)$.
Suppose to the contrary that such intersection occurs. By equating the respective $y$-coordinates with regard to the last two relations from \eqref{inters.that}, we see that this intersection may hold only at symmetric times $t=\pm t_0 \in (-\tau/2, \tau/2)$. Then the respective $x$-coordinates are different since $x_{\circlearrowleft}(\cdot)$ is odd by \eqref{x.even}. This contradiction completes the proof.
\epf
\par
The last claim of Lemma~\ref{lem.typical} entails that the velocity orientation angle $\theta$ arrives at values of the form $\varphi + \psi + 2\pi k$, with $\psi=\pi, \pi/2, -\pi/2$ and an integer $k$, every $\tau$ units of time. At the respective time instants, the robot moves either in the anti-gradient direction or perpendicular to it; see Fig.~\ref{fig.converge}(a).
Thus unrealistic request $\nu > \ov{v}n$ results in motion that involves systematic sideward and backward maneuvers and so is rather ineffective, despite overall gradient climbing.
\par
By the foregoing, the parameter $\mu>0$ of the controller \eqref{c.a} influences neither convergence nor behavioral pattern.
To illuminate the role of $\mu$, we focus on the case of a realistic request $\nu <\ov{v}n$, fix initial location $\bldr_{\text{in}}$, and examine various $\mu>0$. Let $\bldr_\mu(\cdot)$ and $\theta_\mu(\cdot)$ correspond to the respective closed-loop trajectory. Since $\theta(t) := \theta_1(t)$ satisfies equation \eqref{eq.chi} of the closed loop system,
where $\mu:=1$, it is easy to see via the change of the variables $\tau := \mu \tau^\prime, t := \mu t^\prime$ that
$\theta(t):=\theta_1(\mu t)$ solves \eqref{eq.chi} with $\mu$ at hand.
Thus $\theta_\mu(t) = \theta_1(\mu t)$.
It follows that increase of $\mu$ proportionally shortens both time and length of the transient to
the beneficial equilibrium $\theta = \theta_\dagger$, where the requested field growth rate is achieved.
Indeed let $t_0$ be the time for which $\theta_1(\cdot)$ reaches the $\varepsilon$-neighborhood of $\theta_\dagger$, i.e., $|\theta(t) - \theta_\dagger| < \varepsilon\; \forall t \geq t_0$. Then
$\theta_\mu(\cdot)=\theta_1(\mu t)$ evidently does this $\mu$ times faster, i.e., satisfies the same inequality for $t \geq t_0^\mu:=\mu^{-1}t_0$. During the transient to that neighborhood, the deviation of $\bldr_\mu(\cdot)$ from the initial location
\begin{multline*}
\max_{t \in [0, t^\mu_0|} \|\bldr_\mu(t) - \bldr_{\text{in}}\| \overset{t:= \mu^{-1} s}{=\!=\!=\!=} \ov{v} \max_{s \in [0, t_0|} \left\| \int_0^{\mu^{-1} s} \vec{e}[\theta_\mu(\tau)]\, d \tau\right\|
\\
\overset{\tau:= \mu^{-1} \zeta}{=\!=\!=\!=} \mu^{-1} \ov{v} \max_{s \in [0, t_0|} \left\| \int_0^{s} \vec{e}[\theta_1(\zeta)]\, d \zeta\right\| = \mu^{-1} \max_{s \in [0, t_0|} \|\bldr_1(s) - \bldr_{\text{in}}\|
\end{multline*}
is $\mu$ times less than that for $\bldr_1(\cdot)$.
Thus increase of $\mu$ shortens transients in both time and space. This is of interest for generic fields with varying gradients. Large enough $\mu$ gives promise to convert transients into so ``fast'' and ``short'' motions that
they are negligible as compared with the gradient change rate. Since after the transient, the fields value ascends, this should ensure successful overall gradient climbing.
\par
Thus we have acquired first evidence that tuning of the controller \eqref{c.a} should follow two guidelines: firstly, the requested field growth rate $\nu$ should be realistic and secondly, the parameter $\mu$ should be large enough.
To flesh out, specify, and justify this, we proceed to theoretical analysis of the closed-loop behavior in a generic smooth dynamic field $D(\cdot)$.

\section{Convergence to the Maximizer in Generic Dynamic Fields}
\setcounter{equation}{0}
\label{sec.maxx}
We start with notations, which are mainly concerned with characteristics of dynamic fields.
\subsection{Basic notations and field characteristics}
\label{subsec.not}
\begin{figure}[h]
\centering
\scalebox{0.3}{\includegraphics{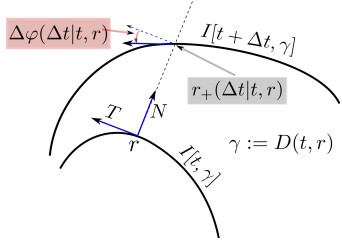}}
\caption{Two close isolines.
}
\label{fig.converge}
\end{figure}
\begin{itemize}
\item $\nabla = \left( \begin{smallmatrix} \frac{\partial}{\partial x} \\ \frac{\partial}{\partial y}  \end{smallmatrix}\right)$ --- the spatial gradient;
    \item $D^{\prime\prime}$ --- the spatial Hessian, i.e., the matrix of the second derivatives with respect to $x$ and $y$;
\item $I(t,\gamma) := \{\bldr: D(t,\bldr) = \gamma \}$ --- the spatial isoline, i.e., the level curve of $D(t,\cdot)$ with the field level $\gamma$;
\item $[T,N]= [T(t,\bldr),N(t,\bldr)]$ --- the (right) Frenet frame of $I[t,\gamma]$ with $\gamma := D(t,\bldr)$ at the point $\bldr$, i.e.,
    $
N = \frac{\nabla D(t,\bldr)}{\|\nabla D (t,\bldr)\|}
$
and $T$ are the unit normal and tangent vectors, respectively;
\item $\varkappa = \varkappa(t,\bldr)$ --- the signed curvature of the spatial isoline;
\item $\bldr_+(\Delta t| t,\bldr)$ --- the nearest (to $\bldr$) point where the ordinate axis of the above Frenet frame of the isoline $I(t,\bldr)$ intersects the isoline with the same field level $\gamma$ observed at time $t=t+\Delta t$; see Fig.~\ref{fig.converge}; 
    \item $p(\Delta t| t,\bldr)$ --- the ordinate of $\bldr_+(\Delta t| t,\bldr)$, i.e., the final normal displacement of the isoline from $t$ to $t+\Delta t$;
\item $\lambda(t,\bldr)$ --- the front velocity of the spatial isoline:
    $
\lambda (t,\bldr):=\lim_{\Delta t\to 0}\frac{p(\Delta t| t,\bldr)}{\Delta t };
$
\item $\alpha(t,\bldr)$ --- the front acceleration of the spatial isoline:
$\alpha(t,r):= \lim_{\Delta t\to 0}\frac{\lambda[t+\Delta t, \bldr_+(\Delta t|t,\bldr)] - \lambda[t,\bldr]}{\Delta t }$;
\item $\Delta \varphi(\Delta t|t,\bldr)$ --- the angular displacement of $T[t+\Delta t, \bldr_+(\Delta t|t,\bldr) ]$ with respect to $T[t,\bldr]$;
\item $\omega(t,\bldr)$ --- the angular velocity of rotation of the spatial isoline:
    $
\omega(t,\bldr):=\displaystyle{\lim_{\Delta t\to 0}}\frac{\Delta \varphi(\Delta t|t,\bldr)}{\Delta t }
$;
\item $\rho(t,\bldr)$ --- the density of isolines:
$
\rho(t,\bldr):=\displaystyle{\lim_{\Delta \gamma \to 0}}\frac{\Delta \gamma}{q(\Delta \gamma| t,\bldr)},
$
where $q(\Delta \gamma| t,\bldr)$ is the ordinate of the nearest to $\bldr$ point where the ordinate axis of the Frenet frame intersects the isoline $I(t| t,\gamma+\Delta \gamma), \gamma:= D(t,\bldr)$; \footnote{This density characterizes the ``number'' of isolines within the unit distance from the basic one $I(t,\gamma)$, where the ``number'' is evaluated by the discrepancy in the values of $D(\cdot)$ observed on these isolines at the given time $t$.}
\item $v_\rho(t,\bldr)$ --- the (proportional) growth rate of the density $\rho$ over time:
\begin{equation}
\label{vrho.def}
v_\rho(t,\bldr):= \frac{1}{\rho(t,\bldr)}\lim_{\Delta t\to 0}\frac{\rho[t+\Delta t, \bldr_+(\Delta t|t,\bldr)]- \rho[t,\bldr]}{\Delta t};
\end{equation}
\item $ \tau_\rho(t,\bldr)$ --- the (proportional) growth rate of the density under a tangential displacement at time $t$:
\begin{equation}
\label{beta_def}
\tau_\rho (t,r):= \frac{1}{\rho(t,\bldr)} \lim_{\Delta s \to 0}\frac{\rho(t, r+ T \Delta s ) -  \rho(t,r)}{\Delta s} ;
\end{equation}
\item $ n_\rho(t,\bldr)$ --- the (proportional) growth rate of the density under a normal displacement at time $t$:
\begin{equation}
\label{njrm_def}
 n_\rho (t,r):= \frac{1}{\rho(t,\bldr)}\lim_{\Delta s \to 0}\frac{ \rho(t, r+ N \Delta s ) -  \rho(t,r)}{\Delta s} ;
\end{equation}
\item $\omega_\nabla(t,\bldr)$ --- the angular velocity of the gradient $\nabla D$ rotation at time $t$ at point $\bldr$.
\end{itemize}
The following lemma explicitly links the above quantities with $D(\cdot)$. Further  $\spr{\cdot}{\cdot}$ is the standard inner product.
\begin{lemma}
\label{lem.relation}
Whenever the field $D(\cdot)$ is twice continuously differentiable in a vicinity of $(t,\bldr)$ and $\nabla D(t,\bldr) \neq 0$, the afore-introduced quantities are well-defined and the following relations hold at $(t,\bldr)$:
 \begin{eqnarray}
 \label{speed1}
 \lambda=-\frac{D'_{t}}{\|\nabla D \|}, \quad   \rho = \|\nabla D\|, \quad  \bldr_+(dt|t,\bldr) = \bldr + \lambda N dt + \so(dt),
 \\
 \label{alpha}
 v_\rho = \frac{\spr{\nabla D^\prime_{t}+\lambda \dd{D}N}{N}}{\|\nabla D \|}, \quad
 \omega= - \frac{\spr{\nabla D^\prime_{t}+\lambda \dd{D}N}{T}}{\|\nabla D\|}, \quad
\alpha=- \frac{D^{\prime\prime}_{tt} + \lambda \left\langle \nabla D^\prime_t ; N \right\rangle}{\|\nabla D\|}  - \lambda  v_\rho ,
\\
 \label{tau}
\varkappa = - \frac{\spr{\dd{D}T}{T}}{\|\nabla D\|}, \quad
\omega_\nabla = - \frac{\spr{\nabla D^\prime_t}{T}}{\|\nabla D\|}
, \quad \tau_\rho =\frac{\spr{D^{\prime\prime}N}{T}}{\left\| \nabla D \right\|} ,
 \quad
 n_\rho =\frac{\spr{D^{\prime\prime}N}{N}}{\left\| \nabla D \right\|} .
\end{eqnarray}
\end{lemma}
\pf
The first claim and \eqref{speed1} are immediate from the implicit function theorem \cite{KrPa02}. We put $r_+(\Delta t):=r_+(\Delta t| t,\bldr)$. Due to \eqref{speed1},
\begin{gather}
\nonumber
\nabla D\left[t+dt,\bldr_+(dt)\right]=
\nabla D[t+dt,\bldr +\lambda N dt+\so(dt)]= \nabla D+[\nabla D^\prime_{t}+\lambda \dd{D}N]dt +\so(dt),
\\
\nonumber
\rho[t+\Delta t, \bldr_+(dt)] =
\left\| \nabla D\left[t+dt,\bldr_+(dt)\right] \right\| = \left\| \nabla D+[\nabla D^\prime_{t}+ \lambda \dd{D} N]dt +\so(dt) \right\|
\\
\nonumber
= \left\| \nabla D\right\| + \frac{\left\langle \nabla D; \nabla D^\prime_{t}+\lambda\dd{D} N \right\rangle}{\left\| \nabla D\right\|} dt + \so(dt)
\\
=
\rho + \left\langle N; \nabla D^\prime_{t}+ \lambda \dd{D} N \right\rangle dt + \so(dt) \;\mid \Rightarrow \text{the first formula in \eqref{alpha}},
\label{expan}
\\
\nonumber
N\left[t+dt,\bldr_+(dt)\right] = \frac{\nabla D\left[t+dt,\bldr_+(dt)\right]}{\left\| \nabla D\left[t+dt,\bldr_+(dt)\right] \right\|} = N + \left[ \frac{\nabla D^\prime_{t}+ \lambda \dd{D} N}{\left\| \nabla D\right\|} - \frac{\nabla D}{\left\| \nabla D\right\|^3} \left\langle \nabla D; \nabla D^\prime_{t}+ \lambda  \dd{D}N \right\rangle  \right] dt+ \so(dt)
\end{gather}
$$
=
N + \frac{1}{\left\| \nabla D\right\|}\left[\nabla D^\prime_{t}+ \lambda \dd{D} N - N \left\langle N; \nabla D^\prime_{t}+ \lambda  \dd{D}N \right\rangle  \right] dt+ \so(dt)
\overset{\text{(a)}}{=} N + \frac{\left\langle T; \nabla D^\prime_{t}+ \lambda  \dd{D}N \right\rangle} {\left\| \nabla D\right\|} T dt   + \so(dt),
$$
where (a) holds since $W = \langle W,T \rangle T + \langle W,N \rangle N \; \forall W \in \real^2$. In the Frenet frame $(T,N)$, we have
\begin{equation*}
N\left[t+dt,\bldr_+(dt)\right] = \left( \begin{array}{c}
-\sin \Delta \varphi(dt|t,\bldr)
\\
\cos \Delta \varphi(dt|t,\bldr)
\end{array}\right)
= N + \left( \begin{array}{c}
-\cos 0
\\
- \sin 0
\end{array}\right) \omega dt + \so(dt)
 =
N - \omega T  dt + \so(dt).
\end{equation*}
By equating the coefficient prefacing $T dt$ in the last two expressions, we arrive at the second formula in \eqref{alpha}. Finally
\begin{multline*}
\lambda(t+dt, \bldr_+(dt)) \overset{\text{\eqref{speed1}}}{=}
-\frac{D^\prime_{t}[t+dt, \bldr_+(dt)]}{\|\nabla D(t+dt, \bldr_+(dt))\|}
\\
\overset{\text{\eqref{expan}}}{=}\lambda - \frac{D^{\prime\prime}_{tt} dt + \left\langle \nabla D^\prime_t ; \bldr_+(dt)-r \right\rangle}{\|\nabla D\|} + D^\prime_t \frac{\left\langle N; \nabla D'_{t}+\lambda D^{\prime\prime} N \right\rangle}{\|\nabla D\|^2} dt +\so(dt)
\\
\overset{\text{\eqref{speed1}}}{=} \lambda - \frac{D^{\prime\prime}_{tt} + \lambda \left\langle \nabla D^\prime_t ; N \right\rangle}{\|\nabla D\|} dt - \lambda  \frac{\left\langle N; \nabla D'_{t}+\lambda D^{\prime\prime} N \right\rangle}{\|\nabla D\|} dt +\so(dt) \overset{\text{(b)}}{=}
\lambda - \frac{D^{\prime\prime}_{tt} + \lambda \left\langle \nabla D^\prime_t ; N \right\rangle}{\|\nabla D\|} dt - \lambda  v_\rho dt +\so(dt),
\end{multline*}
where (b) follows from the first formula in \eqref{alpha}. The definition of $\alpha$ completes the proof of \eqref{alpha}. The first two equations in \eqref{tau} are well known, the third one follows from the transformation
$$
\rho(t, r+ T  ds ) = \left\| \nabla D[t, r + T ds] \right\| = \rho(t,r) + \frac{\spr{D^{\prime\prime}T}{\nabla D}}{\left\| \nabla D \right\|} ds + \so(ds)
= \rho(t,r) + \spr{D^{\prime\prime}T}{N} ds + \so(ds),
$$
the fourth equation in \eqref{tau} is established likewise. \epf
\par
It follows from Lemma~\ref{lem.relation} that $\omega = \omega_\nabla - \lambda \tau_\rho, v_\rho = - \omega_\nabla + \lambda n_\rho$.

\subsection{Assumptions Underlying Theoretical Analysis}\label{subsec.ass}
In the most general setting, the problem at hand comprises problems of
global optimization. In the presence of local
extrema, NP-hardness, this mathematical seal for intractability, was
established for even the simplest classes of such problems
\cite{HorPar95}. Since we do not mean to cope with the challenge of NP-hardness, it is permissible to examine the closed-loop behavior in a field with no local extrema. So our first assumption considers a smooth field with a single global spatial maximizer $\bldr^0(t)$ and no local extrema, which field converges to a finite limit $\gamma_{\inf}(t)$ as $\|r\| \to \infty$.
\begin{assumption}
\label{ass.0a}
The field $D(\cdot)$ is $C^2$-smooth, and there exist continuous functions $\bldr^0(t)$ and $\gamma_{\inf}(t)<D[t,\bldr^0(t)]$ such that $\nabla D(t,\bldr) \neq 0 \; \forall \bldr \neq \bldr^0(t)$ and $D(t,\bldr)\to \gamma_{\inf}(t)$ as $\|\bldr\| \to \infty$.
\end{assumption}
Hence $\max_{\bldr \in \br^2} D(t,\bldr)$ is attained at $\bldr=\bldr^0(t)$, and for given $t$, the field value $D(t,\bldr)$ ranges from $\gamma_{\inf}(t)$ to $D[t,\bldr^0(t)]$.
Assumption~\ref{ass.0a} will be relaxed in Remark~\ref{rem.th}.
\par
For steady fields, no other assumptions are imposed. So the remainder of this subsection addresses only unsteady fields.
\par
If during the search, its conditions, characterized by the field parameters, vary over a very wide, up to infinite range, regular closed-loop behavior can hardly be ensured by a common set of controller parameters. So we will try to inclose the path of the robot in a region $Z_{\text{reg}}$ within which the field parameters from \eqref{speed1}---\eqref{tau} are bounded over space and time. This motivates us to discard locations that are excessively both distant and close to $\bldr^0(t)$, with the latter being unwelcome since, e.g., the isolines typically become highly contorted $\varkappa \approx \infty$ when approaching the maximizer $\bldr^0(t)$.
To reduce technicalities, we pick the region $Z_{\text{reg}}$ so that it is upper and lower bounded by some isolines $I[t,\gamma_-(t)]$ and $I[t,\gamma_+(t)]$ whose levels $\gamma_+(t)> \gamma_-(t), \gamma_\pm(t)\in \big( \gamma_{\inf}(t), D[t,\bldr^0(t)] \big) $ may, however, evolve over time.
The role of $Z_{\text{reg}}$ is to accommodate the transition from $\bldr_{\text{in}}$ to the desired $R_\star$-neighborhood $V_\star(t)$ of the maximizer. So we pick $\gamma_\pm(\cdot)$ so that the region links $\bldr_{\text{in}}$ and $V_\star(t)$. Summarizing, we arrive at the following choice
\begin{multline}
\label{def.zreg}
Z_{\text{reg}}:= \{(t,\bldr) : \gamma_-(t) \leq D(t,\bldr) \leq \gamma_+(t) \}, \qquad \gamma_-(0) < D[0,\bldr_{\text{\rm in}}], \\
\gamma_{\inf}(t)<\gamma_-(t) < \gamma_\star(t) < \gamma_+(t)<D[t,\bldr^0(t)],
\end{multline}
and an intermediate level $\gamma_\star(t)$ is picked so close to $\bldr^0(t)$ that
\begin{equation}
\label{desv.star}
D(t,\bldr) \geq \gamma_\star(t) \Rightarrow \bldr \in V_\star(t).
\end{equation}
We recall that $V_\star(t)$ is the desired vicinity of the maximizer where the robot should be ultimately driven.
\par
For steady fields, the levels $\gamma_i(\cdot)$ are chosen constant, $Z_{\text{reg}}$ has the form $C \times [0,\infty)$, with $C$ compact, and so all field parameters from \eqref{speed1}---\eqref{tau} are bounded over $Z_{\text{reg}}$; their upper bounds will be of interest for controller tuning. For time varying fields, they are bounded over any finite time interval by the same argument. The next assumption is non-void only for dynamic fields and claims that boundedness does not degenerate as $t \to \infty$, which is needed only because the subsequent analysis is asymptotic and so assumes that $t$ does go to $\infty$.
\begin{assumption}
\label{ass.last}
There exists constants  $b_\omega^\nabla$, $b_{\aleph}$ for $\aleph = \rho, \lambda, \omega, \varkappa, v, \alpha, n, \tau$ and $\gamma^0_+, \ov{\gamma}$ such that
\begin{multline}
\label{estim}
\begin{array}{l}
 |\rho| \leq b_\rho, \quad
|\lambda| \leq b_\lambda, \quad |\omega| \leq b_\omega,
\quad |\omega_\nabla| \leq b_\omega^\nabla,
\\
|\varkappa| \leq b_\varkappa, \quad |v_\rho| \leq b_v, \quad |\alpha| \leq b_\alpha, \quad |\tau_\rho| \leq b_\tau, \quad |n_\rho| \leq b_n
\end{array}
\qquad \forall (t,\bldr) \in Z_{\text{reg}}; \\ \gamma_+(t) \leq \gamma^0_+, \quad |\dot{\gamma}_i(t)| \leq  \ov{\gamma} \qquad \forall t, i=\pm,\star.
\end{multline}
\end{assumption}
\par
To advance to the maximizer, the robot should transverse isolines. So it is natural to assume that its speed exceeds that of isolines $\ov{v} > |\lambda|$. A particular isoline with the time-varying level $\gamma=\gamma_\star(t)$ characterizes the desired vicinity of $\bldr^0(t)$ within which the robot should ultimately be kept. So it is also reasonable to assume that the mobility of the robot exceeds that of this isoline. Since the front speed of the latter is $\lambda+\rho^{-1}\dot{\gamma}_\star$, which can be shown similarly to \eqref{speed1}, we assume that $\ov{v}>|\lambda+\rho^{-1}\dot{\gamma}_\star|$. Finally we enhance the both assumptions by unifying them into the requirement $\ov{v} > |\lambda| + \rho^{-1} \ov{\gamma}$ imposed in $Z_{\text{reg}}$. It is always satisfied for steady fields since then $\lambda=0$ and $\ov{\gamma}=0$.
\par
Finally, all assumed strict inequalities are protected from degradation as $t \to \infty$:
\begin{equation}
\label{esstt.grad}
\rho(t,\bldr) \varogreaterthan 0 , \quad v \varogreaterthan |\lambda(t,\bldr)| + \rho^{-1} \ov{\gamma} \; \quad \text{\rm in}\; Z_{\text{reg}}; \qquad \gamma_-(t) \varolessthan \gamma_\star(t) \varolessthan \gamma_+(t) \quad \text{\rm in}\; [0,\infty).
 \end{equation}
Here the notation {\it $f \varogreaterthan g$ in $Z$} is used to express that there exists $\ve >0$ such that $f(t,\bldr) \geq g(t,\bldr) + \ve$ for all $(t,\bldr) \in Z$; the relation $\varolessthan$ is defined likewise. In \eqref{esstt.grad}, the first inequality protects the requirement $\rho(t,\bldr) = \|\nabla D(t,\bldr)\| > 0$ from Assumption~\ref{ass.0a}, whereas the last protects the inequalities $\gamma_-(t) < \gamma_\star(t) < \gamma_+(t)$ from \eqref{def.zreg}.
\par
For steady fields, Assumption~\ref{ass.last} and \eqref{esstt.grad} are automatically fulfilled. It follows from  Lemma~\ref{lem.relation} that for generic unsteady fields satisfying Assumption~\ref{ass.0a} and \eqref{esstt.grad}, boundedness of the field parameters requested in the first row from \eqref{estim} holds whenever the first and second derivatives of the field are bounded on $Z_{\text{reg}}$.
\subsection{Main Theoretical Results}
\label{subsec.mr}
Our first theorem provides an evidence that the proposed control paradigm \eqref{c.a} is intrinsically capable of extremum seeking in an understandable context introduced in the previous subsection.
\begin{theorem}
\label{th.m0} Suppose that \eqref{esstt.grad} and Assumptions~{\rm \ref{ass.0a}, \ref{ass.last}}   are valid. Then there exist parameters $\nu, \mu, d_\ast$ of the controller \eqref{c.a} such that the following claim holds:
\begin{enumerate}[{\rm (i)}]
\item The controller \eqref{c.a} brings the robot to the desired vicinity $V_\star(t)$ of a maximizer in a finite time $t_0$ and keeps it there afterwards: $\bldr(t) \in V_\star(t)$ for $t\geq t_0$.
\end{enumerate}
 Moreover, for any bounded and closed domain $D$ lying in the interior of the set $\{ \bldr: (0,\bldr) \in  Z_{\text{reg}} \}$, there exist common parameters $\nu, \mu, d_\ast$ for which {\rm (i)} holds whenever the initial location $\bldr_{\text{\rm in}} \in D$.
\end{theorem}
The proofs of all claims from this section are given in Section~\ref{sec.proofth}. Now we proceed to discussion of controller tuning.
 Its parameters $d_\ast$ and $\nu$ are chosen prior to $\mu$. Whereas $d_\ast$ is arbitrary, $\nu$ is chosen so that
\begin{equation}
\label{v.astchoice}
\ov{v} \varogreaterthan |\lambda| + \rho^{-1} \nu  \quad  \text{in} \quad Z_{\text{reg}}, \qquad \nu > \ov{\gamma}.
\end{equation}
This is feasible by \eqref{esstt.grad}. For steady fields, $\lambda=0, \ov{\gamma}=0$ and so \eqref{v.astchoice} takes the form $\ov{v}\|\nabla D(\bldr)\| \varogreaterthan \nu >0\; \forall (t,\bldr) \in Z_{\text{reg}}$ by \eqref{speed1}. This conforms to the first guideline stated at the end of Sect.~\ref{sec.linear}: the requested field growth rate $\nu$ should be realistic, i.e., lie in the interval $\dot{d} = \spr{\vec{v}}{\nabla D} \in [-\ov{v}\|\nabla D(\bldr)\|, \ov{v}\|\nabla D(\bldr)\|]$, where the equation holds by \eqref{1}. It can be shown that the first relation from \eqref{v.astchoice} has the same meaning for dynamic fields as well.
By the second and last inequalities in \eqref{v.astchoice} and \eqref{estim}, respectively, the requested field growth rate $\nu$ exceeds the growth rate $\dot{\gamma}_\star(t)$ of the isoline level $\gamma_\star(t)$ that is associated with the desired vicinity of the maximizer.
\par
By the second of the above guidelines, $\mu$ should large enough. To specify this,
we pick lower estimates $\Delta_\nabla, \Delta_\gamma >0$ of the discrepancies in the first, third, and fourth inequalities from \eqref{esstt.grad} and the first strict inequality from \eqref{def.zreg}:
\begin{equation}
\label{rho.lower}
\rho(t,\bldr) \geq \Delta_\nabla  \quad \forall (t,\bldr) \in Z_{\text{reg}}, \quad \gamma_-(t)+\Delta_\gamma \leq \gamma_\star(t) \leq \gamma_+(t) - \Delta_\gamma, \quad D[0,\bldr_{\text{\rm in}}] \geq \gamma_-(0) + \Delta_\gamma/2.
\end{equation}
After invoking the quantities from \eqref{estim}, the choice of
$\mu$ is specified as follows:
\begin{eqnarray}
\label{mu0.choice}
\mu \varogreaterthan \frac{1}{\nu - \ov{\gamma}} \left[-2 \omega - \varkappa v_T + 2 \frac{\tau_\rho \ov{\gamma}}{\rho}  + 2 \frac{v_\rho \ov{\gamma}}{v_T \rho} - \frac{\alpha}{v_T} +  \frac{n_\rho \ov{\gamma}^2}{v_T \rho^2} \right], \quad \text{where} \; v_T := \pm \sqrt{\ov{v}^2 -\left[ \lambda  + \rho^{-1} \ov{\gamma} \right]^2}, \quad \text{in} \quad Z_{\text{reg}};
\\
\label{mu1.choice}
\mu \varogreaterthan \frac{\omega + \tau_\rho (\lambda - \ov{v})}{\rho ( \ov{v} - \lambda) - \nu}
\quad \text{in} \quad Z_{\text{reg}};
\\
\label{mu2.choice}
\mu > \frac{1}{\nu - \ov{\gamma}} \min_{k=1,2,\ldots} \max\{a_1(k), a_2(k)\}, \quad \text{where} \;
\left\{
\begin{array}{l}
a_1(k):=2 b^\nabla_\omega \left[ 1+ \frac{1}{k}\right] + 2 \ov{v} \sqrt{b_\varkappa^2+b_\tau^2} \left[ 2+ \frac{1}{k}\right],
 \\
 a_2(k):= 2b_\rho \frac{(2k+1)\ov{v}+ 2 \pi (k+1) \left( b_\lambda + \Delta_\nabla^{-1} \ov{\gamma}\right)}{\Delta_\gamma}
\end{array}
\right. .
\end{eqnarray}
Since $a_1(x)$ and $a_2(x)$ are descending and ascending functions of $x>0$, respectively, $\min_{k}$ in \eqref{mu2.choice} is attained at the integer either ceiling $\lceil x_{\min} \rceil$ or floor $\lfloor x_{\min} \rfloor$ of the positive root $x_{\min}$ of the quadratic equation $a_1(x) = a_2(x)$.
\par
Now we show that under the recommended choice of the parameters, the control objective is achieved.
\begin{theorem}
\label{th.m1}
Suppose that
\eqref{esstt.grad} and Assumptions~{\rm \ref{ass.0a}, \ref{ass.last}} are valid,
and the controller parameters satisfy \eqref{v.astchoice} and \eqref{mu0.choice}---\eqref{mu2.choice}. Then {\rm (i)} from Theorem~{\rm \ref{th.m0}} is true.
\end{theorem}
\par
Choice of $\nu, \mu$ satisfying \eqref{v.astchoice}, \eqref{mu0.choice}, \eqref{mu1.choice} may proceed from estimation of the respective right-hand sides, based on estimates of the field parameters. For example, \eqref{estim} and \eqref{rho.lower} ensure
\eqref{v.astchoice}, \eqref{mu0.choice}, \eqref{mu1.choice} whenever
\begin{eqnarray*}
\ov{v} >b_\lambda + \frac{\nu}{\Delta_\nabla},
\qquad \mu > \frac{b_\omega + b_\tau (b_\lambda - \ov{v})}{\Delta_\nabla ( \ov{v} - b_\lambda) - \nu} ,
\\
\mu > \frac{1}{\nu - \ov{\gamma}}\left[2 b_\omega + b_{\varkappa} \ov{v} + 2 \frac{b_\tau \ov{\gamma}}{\Delta_\nabla}  + 2 \frac{b_v \ov{\gamma}}{\ov{v}_T \Delta_\nabla} + \frac{b_\alpha}{\ov{v}_T} +  \frac{b_n \ov{\gamma}^2}{\ov{v}_T \Delta_\nabla^2} \right]
, \quad \text{where} \; \ov{v}_T := \pm \sqrt{\ov{v}^2 -\left[ b_\lambda  + \frac{\ov{\gamma}}{\Delta_\nabla} \right]^2}.
\end{eqnarray*}
\par
If the isoline levels are constant, $\ov{\gamma}=0$ and so all small enough $\nu$ and large enough $\mu$ obey the bounds \eqref{v.astchoice}, \eqref{mu0.choice}---\eqref{mu2.choice} from Theorem~\ref{th.m1}. This may be used as a guideline for experimental controller tuning.
\begin{figure}[h]
\centering
\scalebox{0.3}{\includegraphics{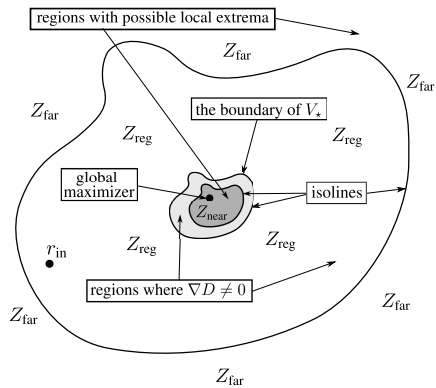}}
\caption{Subdomains of the workspace.
}
\label{fig.sundom}
\end{figure}
\begin{remark}
\label{rem.th}
\rm
Theorems~\ref{th.m0}, \ref{th.m1} remain true if \eqref{esstt.grad} and Assumption~\ref{ass.last} are still stipulated, whereas Assumption~\ref{ass.0a} is relaxed so that the field may have local extrema and be non-smooth in a vicinity $Z_{\text{near}}$ of the global maximizer and in an outlying region $Z_{\text{far}}$; see Fig.~\ref{fig.sundom}. Instead of Assumption~\ref{ass.0a}, it suffices to assume that $D(\cdot)$ is continuous and there exist smooth functions $\gamma_-(t) < \gamma_\star(t) < \gamma_+(t) \in \br$ of time $t$ such that the following statements hold:
\begin{enumerate}[i)]
\item On the set $Z_{\text{reg}}:= \{(t,\bldr) : \gamma_-(t) \leq D(t,\bldr) \leq \gamma_+(t) \}$, the field distribution $D(\cdot)$ is identical to a $C^2$-smooth function defined on a larger and open set, and has no critical points $\nabla D(t,\bldr) \neq 0$;
    \item The set $\{\bldr : \gamma_-(t) \leq D(t,\bldr) \}$ is bounded for any $t$ and $D[0,\bldr_{\text{\rm in}}] > \gamma_-(0)$ ;
\item The implication \eqref{desv.star} is true.
\end{enumerate}
Moreover these theorems remain true even if the regions $Z_{\text{near}} := \{(t,\bldr): D(t,\bldr) > \gamma_{+}(t)\}$ and $Z_{\text{far}} := \{(t,\bldr): D(t,\bldr) < \gamma_{-}(t)\}$ contain point-wise singularities where the field is not defined.
\end{remark}
\par
For continuous fields, $\max_{\bldr \in \br^2} D(t,\bldr)$ is attained and all global maximizers lie in $V_\star(t)$ due to i)---iii).
Thanks to the absence of critical points, the region $Z_{\text{reg}}$ contains no local extrema, whereas they may lie in $Z_{\text{near}}$ and $Z_{\text{far}}$.
\section{Illustrative Examples}
\label{sec.exampl}
\setcounter{equation}{0}
Now we complement the previous discussion in simple yet instructive cases of special and not completely known moving fields. We show that any of these fields satisfies the assumptions of Theorems~\ref{th.m0} and \ref{th.m1}, and
provide an explicit description of controller parameters $(\nu,\mu)$ that guarantee convergence to the field maximizer.
\subsection{Radial Field with a Moving Center}
We consider a scalar field of the form $D(t,\bldr) = c f(\|\bldr-\bldr^0(t)\|)$, where $c>0$ and $\bldr^0(\cdot)$ are unknown, whereas the twice continuously differentiable function $f: [0,\infty) \to \br$ is known,  strictly decaying and convex $f^\prime(z) <0, f^{\prime\prime}(z) \geq 0 \; \forall z>0$. An example of such field is minus the distance to an unknowingly moving target $D(t,\bldr) = - \|\bldr-\bldr^0(t)\|$, in which case $c=1,f(z) = -z$.
Other examples are concerned with some fields caused by emanation of a certain substance from an unknowingly moving point in a homogeneous unbounded medium. If the speed of this point is sufficiently small so that the process is well approximated by the time succession of the instantaneous equilibria, the mathematical model of the field may be given by the above formula, where the unknown $c$ characterizes the emanation rate.
\par
The objective is to display the unknown center $\bldr^0(t)$  of the field by bringing the robot to its $R_\star$-neighborhood $\{\bldr: \|\bldr- \bldr^0(t)\|\leq R_\star\}$ on the basis of the following known estimates
\begin{equation}
\label{known}
 \left\|\bldr_{\text{in}} - \bldr^0 (0)\right\| \leq R_{\text{in}}, \quad \|\dot{\bldr}^0(t)\| \leq v_0, \quad \|\ddot{\bldr}^0(t)\| \leq a_0 \quad \forall t \geq 0, \quad  0 < c_- \leq c \leq c_+ .
\end{equation}
We assume that $R_{\text{in}} > R_\star$ to simplify the formulas.
\begin{proposition}
\label{prop.radial}
Let the robot be faster than the field center $\ov{v}>v_0$, \eqref{known} hold,
and the parameters $\nu,\mu$ of the controller \eqref{c.a} be chosen so that for some $R_- > R_{\text{in}}, 0< R_+ < R_\star$, the following relations are valid:
\begin{eqnarray}
\label{vsats}
0< \nu < c_- |f^\prime(R_-)| (\ov{v}-v_0),
\\
\label{trutru}
\mu   > \nu^{-1} \left[ \frac{\frac{v_0^2}{R_+} + a_0}{\sqrt{\ov{v}^2-v_0^2}} + \frac{\ov{v}+2v_0}{R_+}\right], \quad \mu > \frac{1}{R_+}
\frac{v_0}{ c_- |f^\prime(R_-)| (\ov{v}-v_0) - \nu},
\\
\label{uzhas}
\mu > \frac{4 v_0 +6\ov{v}}{\nu  R_+} , \quad
\mu  > 2 \frac{c_+}{\nu} |f^\prime(R_-)| \frac{3\ov{v}+ 4 \pi  v_0}{\min\{f(R_\star) - f(R_-); f(R_+) - f(R_\star); 2[f(R_{\text{in}}) - f(R_-)]\}}.
\end{eqnarray}
Then the controller \eqref{c.a} brings the robot to the desired $R_\star$-neighborhood of the center $\bldr^0(t)$ in a finite time $t_0$ and keeps it there afterwards $\|\bldr(t) - \bldr^0(t)\| \leq R_\star\; \forall t \geq t_0$.
\end{proposition}
 The described choice of $\nu$ and $\mu$ is possible since in \eqref{vsats}---\eqref{uzhas}, all right-hand sides are finite and positive.
\par
{\bf Proof of Proposition~\ref{prop.radial}:} We take $\gamma_i(t) \equiv cf(R_i), i=\star, \pm$ in \eqref{def.zreg} and invoke that $\Phi_{\frac{\pi}{2}} N$ is the counter-clockwise rotation through the angle $\pi/2$. An elementary calculus exercise based on Lemma~\ref{lem.relation} shows that
\begin{gather*}
N = -\frac{\bldr - \bldr^0}{\|\bldr - \bldr^0\|}, \quad T = - \Phi_{\frac{\pi}{2}} N , \quad
\lambda = \spr{\dot{\bldr}^0}{N}, \quad \rho = c |f^\prime|,  \quad \varkappa =  \frac{1}{\|\bldr-\bldr^0\|}, \quad \tau_\rho =0, \quad n_\rho = \frac{f^\pp}{|f^\prime|},
\\
v_\rho = 0,
\quad \omega = \omega_\nabla = - \frac{\spr{\dot{\bldr}^0}{T}}{\|\bldr-\bldr^0\|},
    \quad \alpha = \spr{\ddot{\bldr}^0}{N} +  \frac{\spr{\dot{\bldr}^0}{T}^2}{\|\bldr - \bldr^0\|}.
\end{gather*}
Now $Z_{\text{reg}} = \{(t,\bldr) : R_+ \leq \|\bldr - \bldr^0(t)\| \leq R_-\}$ is a time-varying annulus and so in \eqref{estim}, \eqref{rho.lower},
\begin{eqnarray*}
\ov{\gamma} =0 , \quad b_\omega = b^\nabla_\omega = \frac{v_0}{R_+}, \quad b_\varkappa = \frac{1}{R_+}, \quad b_\tau =0, \quad b_\rho = c_+ |f^\prime(R_+)|, \quad b_\lambda = v_0,
\\
\Delta_\nabla = c_- |f^\prime(R_-)|, \quad
\Delta_\gamma = c_- \min\{f(R_\star) - f(R_-); f(R_+) - f(R_\star); 2[f(R_{\text{in}}) - f(R_-)]\}.
\end{eqnarray*}
 Hence the assumptions of Theorems~{\rm \ref{th.m0}} and {\rm \ref{th.m1}} are fulfilled in the form given by Remark~\ref{rem.th}, whereas  \eqref{v.astchoice} holds by \eqref{vsats}.
Inequalities \eqref{mu0.choice} and \eqref{mu1.choice} take the form
\begin{equation*}
\mu \nu  \varogreaterthan \left| \frac{ \frac{\spr{\dot{\bldr}^0}{T}^2}{\|\bldr - \bldr^0\|} + \spr{\ddot{\bldr}^0}{N}}{\sqrt{\ov{v}^2 - \lambda^2}} + \frac{ \sqrt{\ov{v}^2 - \lambda^2}}{\|\bldr-\bldr^0\|} \right| + 2 \frac{\spr{\dot{\bldr}^0}{T}}{\|\bldr-\bldr^0\|}, \quad \mu \varogreaterthan - \frac{\spr{\dot{\bldr}^0}{T}}{\|\bldr-\bldr^0\|} \frac{1}{\rho ( \ov{v} - \lambda) - \nu} \quad \text{in} \; Z_{\text{reg}}
\end{equation*}
and follow from \eqref{trutru} since $\spr{\dot{\bldr}^0}{T}$ takes values in $[-\sqrt{v_0^2-\lambda^2}, \sqrt{v_0^2-\lambda^2}]$ for given $\lambda$.
Relations \eqref{mu2.choice} with $k=1$ are implied by \eqref{uzhas}. Theorem~\ref{th.m1} and Remark~\ref{rem.th} complete the proof. \epf
\par
Now we discuss in more details a particular case of the problem at hand.
\subsection{Locating and Escorting a Moving Target Based on Range-Only Measurements}
In this case, $D[t,\bldr] = - \|\bldr - \bldr^0(t)\|, c=c_-=c_+ =1, f(z)=-z$, and the field center $\bldr^0(t)$ is at the target.
\begin{proposition}
\label{prop.target}
Let the target $\bldr^0(t)$ and the robot initial location $\bldr_{\text{in}}$ satisfy the first three inequalities from \eqref{known},
and the robot be faster than the target $\ov{v}>v_0$. Then the conclusion of Proposition~{\rm \ref{prop.radial}} remains true if
\begin{multline}
\label{mov.targ}
0< \nu < \ov{v}-v_0,
\\ \mu > \max\left\{
\frac{2v_0}{ R_\star (\ov{v}-v_0 - \nu)};  \frac{12 \ov{v}+ 8   v_0}{\nu R_\star}; \frac{16 \pi v_0 +12\ov{v}}{ \nu R_\star}; \frac{a_0 }{\nu \sqrt{\ov{v}^2-v_0^2}} + \frac{2}{R_\star \nu} \left[\ov{v}+2v_0 + \frac{v_0^2}{\sqrt{\ov{v}^2-v_0^2}} \right] \right\}.
\end{multline}
\end{proposition}
\pf Now \eqref{vsats}---\eqref{uzhas} take the form
\begin{eqnarray*}
0< \nu < \ov{v}-v_0,
\quad
\mu  \nu >    \frac{a_0 }{\sqrt{\ov{v}^2-v_0^2}} + \frac{1}{R_+} \left[\ov{v}+2v_0 + \frac{v_0^2}{\sqrt{\ov{v}^2-v_0^2}} \right], \quad \mu > \frac{1}{R_+}
\frac{v_0}{ \ov{v}-v_0 - \nu},
\\
\mu \nu > \frac{4 v_0 +6\ov{v}}{R_+} , \quad
\mu \nu > 2 \frac{3\ov{v}+ 4 \pi  v_0}{\min\{R_--R_\star; R_\star- R_+; 2(R_--R_{\text{in}})\}}.
\end{eqnarray*}
The freedom to arbitrarily pick $R_-$ permits us to let $R_- \to \infty$ here, which results in
$$
0< \nu < \ov{v}-v_0,
\quad
\mu  \nu >   \frac{a_0 + \frac{v_0^2}{R_+}}{\sqrt{\ov{v}^2-v_0^2}} + \frac{\ov{v}+2v_0}{R_+}, \quad \mu > \max \left\{ \frac{a}{R_+} ; \frac{b}{R_\star- R_+} \right\}.
$$
where $a:= \max\left\{
\frac{v_0}{ \ov{v}-v_0 - \nu}; \frac{4 v_0 +6\ov{v}}{\nu} \right\}, b:= 2 \frac{3\ov{v}+ 4 \pi  v_0}{\nu} $. By picking $R_+:= R_\star/2$, we arrive at  \eqref{mov.targ}.
Proposition~\ref{prop.radial} completes the proof. \epf
\section{Proofs of Theorems~\ref{th.m0}, \ref{th.m1} and Remark~\ref{rem.th}}
\label{sec.proofth}
\setcounter{equation}{0}
We start with a purely technical exercise in calculus, which addresses quantities introduced in subsect.~\ref{subsec.not}.
\begin{lemma}
The following relations hold:
\begin{gather}
\label{lrho}
\lambda (t,r+T ds) = \lambda + \omega ds + \so (ds), \quad \lambda (t,r+N ds) = \lambda - v_\rho ds + \so (ds);
\\
\label{frenet-serrat}
N[t,r+T ds] = N - \kappa T ds +\so(ds), \quad T[t,r+T ds] = T + \kappa N ds +\so(ds);
\\
\label{nnn}
N(t,r+N ds)= N + \tau_\rho T ds + \so(ds), \quad T(t,r+N ds)= T - \tau_\rho N ds + \so(ds),
\\
\label{n=omega}
N[t+dt, \bldr_+(dt)] = N - \omega T dt +\so(dt), \quad T[t+dt, \bldr_+(dt)] = T + \omega N dt +\so(dt),
\end{gather}
where $\bldr_+(\Delta t):=\bldr_+(\Delta t|t,\bldr)$.
\end{lemma}
\pf
Formulas \eqref{lrho} are justified by the following observations:
\begin{multline*}
\lambda(t,r+T ds) \overset{\text{\eqref{speed1}}}{=} -\frac{D^\prime_{t}(t,r+T ds)}{\|\nabla D(t, r+ T ds)\|} = \lambda(t,r) - \frac{\spr{\nabla D^\prime_{t}}{T}}{\|\nabla D\|} ds + D^\prime_t \frac{\spr{D^{\prime\prime}N}{T} }{\|\nabla D\|^2} ds + \so(ds)
\overset{\text{\eqref{speed1}}}{=} \lambda(t,r) - \frac{\spr{\nabla D^\prime_{t}}{T}}{\|\nabla D\|} ds
\end{multline*}
\begin{multline*}
- \lambda \frac{\spr{D^{\prime\prime}N}{T} }{\|\nabla D\|} ds + \so(ds)
=
\lambda(t,r) - \frac{\spr{\nabla D^\prime_{t}+ \lambda D^{\prime\prime}N}{T}}{\|\nabla D\|} ds  + \so(ds)
\overset{\text{\eqref{alpha}}}{=} \lambda + \omega ds + \so (ds);
\\
\lambda(t,r+N ds) = \lambda(t,r) - \frac{\spr{\nabla D^\prime_{t}+ \lambda D^{\prime\prime}N}{N}}{\|\nabla D\|} ds + \so(ds) \overset{\text{\eqref{alpha}}}{=}
\lambda - v_\rho ds + \so (ds).
\end{multline*}
Formulas \eqref{frenet-serrat} are the classic Frenet-Serrat equations. Furthermore
\begin{multline*}
N(t,r+N ds) = \frac{\nabla D(t,r+N ds)}{\|\nabla D(t,r+N ds)\|} = N + \frac{D^\pp N}{\|\nabla D\|}ds - \nabla D \frac{\spr{D^\pp N}{\nabla D}}{\|\nabla D\|^3}ds + \so(ds)
 \\
 = N + \frac{D^\pp N - N \spr{D^\pp N}{N} }{\|\nabla D\|}ds + \so(ds)
 = N + \frac{\spr{D^\pp N}{T} }{\|\nabla D\|}T ds + \so(ds) \overset{\text{\eqref{tau}}}{=} N + \tau_\rho T ds + \so(ds),
  \end{multline*}
  which yields the first formula in \eqref{nnn}. This also implies the second formula since $N = \Phi_{\frac{\pi}{2}} T, T = - \Phi_{\frac{\pi}{2}} N$. Formulas \eqref{n=omega} follow from the definition of $\omega$. \epf
\par
From now on, we examine the robot driven by the control law
\eqref{c.a}. Assumption~\ref{ass.0a} is taken in the generalized form, i.e., as i)---iii) from Remark~\ref{rem.th}.
Prior to analysis, we note that thanks to \eqref{v.astchoice}, \eqref{mu0.choice}, along with \eqref{estim}, \eqref{esstt.grad}, there exist
\begin{equation}
\label{delta.choice}
\delta \in (0,\nu-\ov{\gamma}), \qquad \delta \approx \nu - \ov{\gamma}
\end{equation}
and $\ve_i>0$ such that everywhere in $Z_{\text{reg}}$,
\begin{gather}
\label{ve.2}
\ov{v} \geq |\lambda| + v_\Delta + \ve_1, \quad \text{where} \quad v_\Delta:= \rho^{-1}(\nu - \delta),
\\
\label{ve.1}
\mu\delta \geq -2 \omega - \varkappa v_T^\prime + 2 \tau_\rho v_\Delta + 2 v_\rho \frac{v_\Delta}{v_T^\prime} - \frac{\alpha}{v_T^\prime} + n_\rho \frac{v_\Delta^2}{v_T^\prime} + \ve_2, \quad \text{where} \; v_T^\prime := \pm \sqrt{\ov{v}^2 -\left[ \lambda  + v_\Delta \right]^2},
\end{gather}
and \eqref{mu2.choice} holds with $\delta$ put in place of $\nu$.
We also observe that the closed-loop system obeys the equations
\begin{equation}
\label{c.a2}
\dot{\bldr} = \ov{v} \vec{e}(\theta), \quad \dot{\theta} =
\mu\big[\dot{d} - \nu \big] , \quad \dot{d} = \ov{v} \spr{\nabla D(t,\bldr)}{\vec{e}(\theta)} + \frac{\partial D}{\partial t}(t,\bldr).
\end{equation}
The next two lemmas display the key features of the control law
\eqref{c.a}. The first lemma shows that under certain conditions, $\dot{d}= \nu-\delta \Rightarrow \ddot{d}>0$. This implies that $\dot{d}$ cannot trespass $\nu-\delta$ from above and so if the inequality $\dot{d} \geq \nu-\delta$ is achieved, it will be kept true afterwards. By \eqref{estim} and \eqref{delta.choice}, this implies that $\dot{d} > |\dot{\gamma}_\star|$ and so $d(t)$ not only ascends but also eventually overtakes $\gamma_\star(t)$. The last means that the robot enters the desired vicinity of the maximizer by \eqref{desv.star}.
\begin{lemma}
\label{lem.basic}
For all $(t,\bldr) \in Z_{\text{reg}}$, the following implication holds
\begin{equation}
\label{impl}
\dot{d}= \nu -\delta \Rightarrow
\begin{cases}
\ddot{d} > 0 & \text{\rm if}\; \sigma:= \spr{T}{\vec{e}(\theta)} <0
\\
\ddot{d} <0 & \text{\rm if}\; \sigma >0
\end{cases}
.
\end{equation}
\end{lemma}
\pf
With start by noting that
\begin{gather}
\nonumber
\nu -\delta = \dot{d} =  D'_{t}+\ov{v}<\nabla D,\vec{e}>
 \overset{\text{\eqref{speed1}}}{=}\rho \left( - \lambda + \ov{v} \left\langle N, \vec{e} \right\rangle \right)  \Rightarrow  \left\langle N, \vec{e} \right\rangle = \frac{\lambda}{\ov{v}} + \frac{v_\Delta}{\ov{v}} , \quad \text{where} \quad v_\Delta:= \rho^{-1}(\nu - \delta),
 \\
 \label{e.expand}
 \Rightarrow \vec{e} = \frac{1}{\ov{v}}\Big[ \left( \lambda   + v_\Delta \right) N + \underbrace{\sgn \sigma \sqrt{\ov{v}^2 -\left[ \lambda  + v_\Delta \right]^2}}_{v_T}  T \Big];
 \\
 \label{rplus}
 r(t +dt) = r(t) + \ov{v} \vec{e} dt +\so(dt) = \underbrace{r(t) + \lambda N \, dt}_{= \bldr_+(dt) +\text{\tiny $\mathcal{O}$}(dt) \text{by \eqref{speed1}}} +  [ v_\Delta N +  v_T  T ] dt +\so(dt);
 \\
  \nonumber
 \dot{d}(t+dt) - \dot{d}(t) = \rho[t+dt,\bldr(t+dt)] \left( - \lambda [t+dt,\bldr(t+dt)] + \ov{v} \left\langle N[t+dt,\bldr(t+dt)], \vec{e}(t+dt) \right\rangle \right)
\end{gather}
 \begin{gather}
 \nonumber
 - \rho[t,\bldr(t)] \left( - \lambda [t,\bldr(t)] + \ov{v} \left\langle N[t,\bldr(t)], \vec{e}(t) \right\rangle \right)
\\
 \nonumber
 = \underbrace{\big\{ \rho[t+dt,\bldr(t+dt)] - \rho[t,\bldr(t)]\big\}}_{a}\underbrace{( - \lambda + \ov{v} \left\langle N, \vec{e} \right\rangle)}_{= v_\Delta}
 - \rho \underbrace{\left\{ \lambda [t+dt,\bldr(t+dt)] - \lambda [t,\bldr(t)]\right\}}_{b}
 \\
 \nonumber
 + \ov{v} \rho  \big\langle \underbrace{N[t+dt,\bldr(t+dt)] - N[t,\bldr(t)]}_{c1}, \vec{e} \big\rangle
 + \ov{v} \rho \big\langle N, \underbrace{\vec{e}(t+dt) - \vec{e}(t)}_{c2} \big\rangle + \so(dt);
 \\
 \nonumber
 a \overset{\text{\eqref{rplus}}}{=} \big\{ \rho[t+dt,\bldr_+(dt) + v_\Delta N dt + v_T T dt + \so(dt)] - \rho[t,\bldr(t)] \overset{\text{\eqref{vrho.def}--\eqref{njrm_def}}}{=\!=\!=\!=\!=\!=}
 \rho [v_\rho + v_\Delta n_\rho + v_T \tau_\rho ]dt + \so(dt).
\end{gather}
Similarly and with regard to the definition of the front acceleration $\alpha$,
\begin{gather}
\nonumber
b \overset{\text{\eqref{lrho}}}{=\!=\!=} [\alpha + \omega v_T - v_\rho v_\Delta] dt + \so(dt), \quad c1 \overset{\text{\eqref{frenet-serrat}--\eqref{n=omega}}}{=\!=\!=\!=\!=\!=} [- \omega + \tau_\rho v_\Delta - \varkappa v_T] T dt + \so(dt),
\\
 \nonumber
   c2 \overset{\text{\eqref{1}}}{=} \dot{\theta} \Phi_{\frac{\pi}{2}} \vec{e} dt +\so(dt) \overset{\text{\eqref{e.expand}}}{=}  \frac{\dot{\theta}}{\ov{v}} \left[- \left( \lambda   + v_\Delta \right) T + v_T  N \right] dt + \so(dt) .
\end{gather}
Summarizing with regard to \eqref{e.expand} and the expansion $\dot{d}(t+dt) - \dot{d}(t) = \ddot{d}(t) dt + \so(dt)$, we see that
\begin{gather}
\nonumber
\ddot{d}(t) =  \rho \big[ (v_\rho + v_\Delta n_\rho + v_T \tau_\rho )v_\Delta - (\alpha + \omega v_T - v_\rho v_\Delta) + (- \omega + \tau_\rho v_\Delta - \varkappa v_T)v_T + v_T  \dot{\theta} \big]
\\
\nonumber
=
\rho v_T \left[ \dot{\theta} - 2 \omega - \varkappa v_T+2 \tau_\rho v_\Delta + 2 v_\rho \frac{v_\Delta}{v_T} - \frac{\alpha}{v_T} + n_\rho \frac{v_\Delta^2}{v_T}\right], \qquad \text{where} \quad \dot{\theta} \overset{\text{\eqref{c.a2}}}{=} \mu (\dot{d} - \nu) \overset{\dot{d}=\nu - \delta}{=\!=\!=\!=\!=}  - \mu \delta .
\end{gather}
 The proof is completed by noting that the expression inside $[\ldots]$ does not exceed $- \ve_1$ by \eqref{ve.1},
and $v_T \overset{\text{\eqref{e.expand}}}{=} \ov{v} \spr{T}{\vec{e}}$. \epf
\par
By Lemma~\ref{lem.basic}, the key inequality $\dot{d} \geq \nu-\delta$ is maintained whenever $\sigma <0$. The next lemma implies that while the inequality is true, $\sigma$ cannot trespass zero from below and so the condition $\sigma <0$ is kept true.
\begin{lemma}
\label{lem.sigma}
For $ \sigma(t) = \spr{T[t,\bldr(t)]}{\vec{e}[\theta(t)]} $, the following implication holds
$$
\sigma =0 \wedge \dot{d} \geq \nu - \delta \Rightarrow \dot{\sigma}<0.
$$
\end{lemma}
\pf
We first note that $\sigma=0 \Rightarrow \vec{e} = \pm N$. If $\vec{e} = - N$, we have
$$
\dot{d} = D^\prime_t+ \ov{v} \spr{\nabla D}{\vec{e}}= D^\prime_t - \ov{v} \rho
\overset{\eqref{speed1}}{=\!=} - \rho (\lambda + \ov{v}) \overset{\text{\eqref{esstt.grad}}}{<} 0,
$$
in violation of $\dot{d} \geq \nu - \delta >0$. Hence $\vec{e} = N$.
Similarly to \eqref{rplus},
\begin{gather}
\nonumber
r(t +dt) = r(t) + \ov{v} \vec{e} dt +\so(dt) = r(t) + \ov{v} N dt +\so(dt)  = r_+(t) + (\ov{v} -\lambda) N dt +\so(dt) ,
\\
\nonumber
T[t+dt,\bldr(t+dt)] - T[t,\bldr(t)] \overset{\text{\eqref{nnn},\eqref{n=omega}}}{=\!=\!=\!=\!=\!=} [\omega + \tau_\rho (\lambda - \ov{v})] N dt + \so(dt),
\\
\nonumber
\vec{e}[\theta(t+dt)] - \vec{e}[\theta(t)] \overset{\text{\eqref{1}}}{=} \dot{\theta} \Phi_{\frac{\pi}{2}} \vec{e} dt +\so(dt) =  \dot{\theta} \Phi_{\frac{\pi}{2}} N dt +\so(dt) = - \dot{\theta} T dt +\so(dt),
\\
\nonumber
\sigma(t+dt) - \sigma(t) =  \spr{T[t+dt,\bldr(t+dt)] - T[t,\bldr(t)]}{\vec{e}}
+  \spr{T}{\vec{e}[\theta(t+dt)] - \vec{e}[\theta(t)]}
\\
\nonumber
=
\spr{[\omega + \tau_\rho (\lambda - \ov{v})] N}{\vec{e}} dt - \spr{T}{\dot{\theta} T} dt + \so(dt)
= \big[ \omega + \tau_\rho (\lambda - \ov{v}) - \dot{\theta} \big]dt + \so(dt);
\\
\label{dott.dd}
\dot{d} = D^\prime_t+ \ov{v} \spr{\nabla D}{\vec{e}} = D^\prime_t + \ov{v} \rho = \rho(\ov{v} - \lambda );
\quad \dot{\theta}  \overset{\eqref{c.a}}{=\!=} \mu (\dot{d} - \nu) = \mu \rho ( \ov{v} - \lambda - \nu/\rho );
\\
\nonumber
\dot{\sigma} = \omega + \tau_\rho (\lambda - \ov{v}) - \mu \rho ( \ov{v} - \lambda - \nu/\rho ) \overset{\text{\eqref{mu1.choice}}}{<}0. \quad \text{\epf}
\end{gather}
\par
Up to minor details, the next claim shows that the robot reaches the desired vicinity of the maximizer provided that initially the afore-mentioned inequalities $\dot{d}(t) \geq \nu - \delta$ and $\sigma(t) \leq 0$ hold.
\begin{cor}
\label{cor.mcase}
Suppose that $\dot{d}(t) \geq \nu - \delta, \sigma(t) \leq 0$, and $[t,\bldr(t)] \in Z_{\text{reg}}$ at $t=t_0$. In a finite time, the robot reaches the curve $I[t,\gamma_+(t)]$, with maintaining the above inequalities and inclusion true.
\end{cor}
\pf
If $\bldr(t_0) \in I[t_0,\gamma_+(t_0)]$, the claim is trivial. Otherwise, $D[t,\bldr(t)] < \gamma_+(t), \dot{d} \geq \nu - \delta >0$ at $t=t_0$, and  $\dot{d}(t_0) = \nu - \delta \Rightarrow \ddot{d}(t_0) >0$ by Lemma~\ref{lem.basic}, whereas $\sigma(t_0)=0 \Rightarrow \dot{\sigma}(t_0) <0$ by Lemma~\ref{lem.sigma}. So the open set
$$
S:= \{t>t_0: \gamma_-(t)< d(t) < \gamma_+(t), \dot{d}(t) > \nu -\delta, \sigma(t) <0 \}
$$
contains an interval of the form $(t_0, t_\ast), t_\ast > t_0$. Let us consider the maximal of such $t_\ast$'s. Here $t_\ast \neq \infty$ since otherwise $\gamma_+(t)> d(t) \geq d(t_0) + (\nu-\delta)(t-t_0) \to \infty$ as $t \to \infty$, in violation of the penultimate inequality from \eqref{estim}. Hence $t_\ast < \infty$ and either (a) $\gamma_-(t_\ast) = d(t_\ast)$ or (b) $\sigma(t_\ast) =0$, or (c) $\sigma(t_\ast)<0$ and $\dot{d}(t_\ast) = \nu - \delta$, or (d) $d(t_\ast) = \gamma_+(t_\ast)$.
However (b) and (c) are impossible by Lemmas~\ref{lem.sigma} and \ref{lem.basic}, respectively, whereas
$$
d(t_0) \geq \gamma_-(t_0) \wedge \dot{d}(t) > \nu-\delta \overset{\text{\eqref{delta.choice}}}{>} \ov{\gamma} \overset{\text{\eqref{estim}}}{\geq} \dot{\gamma}_-(t) \Rightarrow d(t_\ast) > \gamma_-(t_\ast)
$$
and so (a) does not hold either. Thus (d) holds, which completes the proof. \epf
\par
Now we show that even if the above inequalities $\dot{d}(t) \geq \nu - \delta$ and $\sigma(t) \leq 0$ do not initially hold, they become true at some time instant. To this end, we need some technical facts concerned with motion in the situation where some of them is violated.
\begin{lemma}
\label{lem.diberg}
Suppose that $\dot{d}(t) \leq \nu - \delta, [t,\bldr(t)] \in Z_{\text{reg}} \; \forall t \in \Delta = [t_0,t_1]$. For $t \in \Delta$, the vector $\vec{e}(\theta)$ rotates clockwise with the angular velocity $ \dot{\theta} \leq - \mu \delta$, and the space deviation from the initial location obeys the inequality
\begin{multline}
\label{ineq.dist}
\|\bldr(t) - \bldr(t_0)\|
\leq q(\varphi):= \frac{2v}{\mu \delta} \left\lfloor \frac{\varphi}{2 \pi} \right\rfloor + \frac{\ov{v}}{\mu \delta} \left[ 1 - \cos \min\{ \Lbag \varphi \Rbag; \pi \}\right], \\ \text{\rm where} \quad \varphi := |\theta(t)-\theta(t_0)|, \quad  \Lbag \varphi \Rbag := \varphi - 2 \pi \left\lfloor \frac{\varphi}{2 \pi} \right\rfloor.
\end{multline}
\end{lemma}
\pf
Without any loss of generality, we assume that $t_0=0, \bldr(0)=0, \theta(0) =0$.
The first claim of the lemma is immediate from the second equation in \eqref{c.a2}:
$
 \dot{\theta} = \mu \left[ \dot{d} - \nu \right] \overset{\dot{d} \leq \nu -\delta}{\leq} - \mu \delta.
$
So $s:= -\theta$ can be taken as a new independent variable:
$$
\frac{d \bldr }{d s} = u \vec{e}(-s), \quad u:= - \frac{\ov{v}}{\dot{\theta}} \in \left[ 0 , \frac{\ov{v}}{\mu\delta} \right].
$$
The squared distance $\|\bldr(\varphi)\|^2$ does not exceed the maximal value of $\mathcal{I}$ in the following optimization problem:
\begin{equation*}
\mathcal{I}:= \|\bldr(\varphi)\|^2 \to \max \quad \text{subject to} \quad \frac{d \bldr }{d s} = u \vec{e}(-s) \; s \in [0,\varphi], \quad \bldr(0)=0, \quad u (s) \in \left[ 0 , \ov{v}(\mu\delta)^{-1} \right].
\end{equation*}
Its solution $\bldr^0(\cdot), u^0(\cdot)$ exists and obeys the Pontryagin's maximum principle \cite{Vin00}: there exists a smooth function $\psi (s) \in \br^2$ such that
\begin{multline*}
u^0(s) = \text{arg\,max}_{u \in [0, \ov{v}/(\mu\delta)]} u \psi^\trs \vec{e}(-s)
= \begin{cases}
 \ov{v}/(\mu\delta) & \text{if}\; \psi^\trs \vec{e}(-s) >0
 \\
 0 & \text{if}\; \psi^\trs \vec{e}(-s) < 0
  \\
 \text{unclear} & \text{if}\; \psi^\trs \vec{e}(-s) = 0
\end{cases}, \\ \text{where} \;
\begin{array}{l}
\frac{d \psi }{ds} = - \frac{\partial}{\partial \bldr} u^0 \psi^\trs \vec{e}(-s) =0 \Rightarrow \psi = \text{const}
\\
\psi = \psi(\varphi) = 2 \bldr^0(\varphi),
\end{array}.
\end{multline*}
If $\bldr^0(\varphi) = 0$, \eqref{ineq.dist} is evident.
Let $\bldr^0(\varphi)\neq 0$. Then $\psi \neq 0$ and
so as $s$ progresses, $u^0(s)$ interchanges the values $0$ and $\ov{v}/(\mu\delta)$, each taken on an interval of length $\pi$ possibly except for  extreme intervals whose lengths do not exceed $\pi$. Inequality \eqref{ineq.dist} results from computation of $\|\bldr(\varphi)\|$ for such controls $u(\cdot)$, along with picking the maximum among these results.
\epf
\par
The next step is a technical study of isolines with time-varying field levels, like those bounding the desired vicinity of the maximizer and the operational zone $Z_{\text{reg}}$ by i)---iii) of Remark~\ref{rem.th}.
For any simply connected domain $\mathfrak{D} \subset Z_{\text{reg}}$, the angle $\beta$ of rotation of the vector-field $\nabla D$ along a curve from $\mathfrak{D}$ is uniquely determined by the ordered pair of the ends of this curve since $\nabla D (t,\bldr) \neq 0 \; \forall (t,\bldr) \in \mathfrak{D}$ by \eqref{esstt.grad}.  Let $\beta_\nabla(\mathfrak{D})$ be the maximum of $|\beta|$ over all pairs in $\mathfrak{D}$.
\begin{lemma}
\label{lemrroott}
Suppose that $\mathfrak{D}:=\Delta \times C \subset Z_{\text{reg}}$, where $\Delta$ is an interval and $C \subset \br^2$ is a convex set. Then
$$
\beta_\nabla(\mathfrak{D}) \leq \sup_{p \in \mathfrak{D}} |\omega_\nabla(p)| \times \text{\bf diam} \Delta + \sup_{p \in \mathfrak{D}} \sqrt{\varkappa^2+\tau_\rho^2} \times \text{\bf diam} C,
$$
where $\text{\bf diam}$ is the diameter of the set, i.e., the supremum of distances between two its elements.
\end{lemma}
\pf
Let $p_i=(t_i, \bldr_i) \in \mathfrak{D}:=\Delta \times C, i=0,1$ and $t_0 \leq t_1$ for the definiteness. In $\mathfrak{D}$, we introduce two parametric curves $\zeta(s):= [s, \bldr_0], s \in [t_0,t_1]$ and $\xi(s):= [t_1, (1-s)\bldr_0 + s \bldr_1], s \in [0,1]$, connecting the point $[t_0,\bldr_0]$ with $[t_0,\bldr_1]$ and $[t_0,\bldr_1]$ with $[t_1,\bldr_1]$, respectively. The concatenation of these curves connects $p_1$ with $p_2$ and lies in $\mathfrak{D}$. So the angle of the gradient rotation when going from $p_1$ to $p_2$ inside $\mathfrak{D}$ is the angle of rotation $\sphericalangle \beta_\zeta $ when going along $\zeta(\cdot)$ plus the angle of rotation $\sphericalangle \beta_\xi $ when going along $\xi(\cdot)$. Let $\varphi[t,\bldr]$ stands for the orientation angle of $\nabla D [t,\bldr]$ in the absolute coordinate frame. Then
\begin{eqnarray*}
\frac{d}{ds} \varphi[\zeta(s)] = - \frac{\spr{\nabla D^\prime_t [\zeta(s)]}{T[\zeta(s)]}}{\|\nabla D[\zeta(s)]\|} \overset{\text{\eqref{tau}}}{=\!=} \omega_\nabla[\zeta(s)] \Mapsto |\sphericalangle \beta_\zeta | = \left| \int_{t_0}^{t_1} \omega_\nabla[\zeta(s)] \;ds \right| \leq \sup_{p \in \mathfrak{D}} |\omega_\nabla(p)| (t_1-t_0);
\\
\frac{d}{ds} \varphi[\xi(s)] = - \frac{\spr{D^\pp [\zeta(s)][\bldr_1-\bldr_0]}{T[\zeta(s)]}}{\|\nabla D[\zeta(s)]\|}
= - \frac{\spr{D^\pp [\zeta(s)]T[\zeta(s)]}{T[\zeta(s)]}}{\|\nabla D[\zeta(s)]\|} \spr{\bldr_1-\bldr_0}{T[\zeta(s)]} \\ - \frac{\spr{D^\pp [\zeta(s)]T[\zeta(s)]}{N[\zeta(s)]}}{\|\nabla D[\zeta(s)]\|} \spr{\bldr_1-\bldr_0}{N[\zeta(s)]} \overset{\text{\eqref{tau}}}{=\!=}
 +\varkappa \spr{\bldr_1-\bldr_0}{T[\zeta(s)]} - \tau_\rho \spr{\bldr_1-\bldr_0}{N[\zeta(s)]} \Mapsto
\\
\left| \frac{d}{ds} \varphi[\xi(s)] \right| \leq \|\bldr_1-\bldr_0\| \sqrt{\varkappa^2+\tau_\rho^2} \Mapsto |\sphericalangle \beta_\xi | \leq
\int_0^1 \left| \frac{d}{ds} \varphi[\xi(s)] \right| \; ds \leq \|\bldr_1-\bldr_0\| \sup_{p \in \mathfrak{D}} \sqrt{\varkappa^2+\tau_\rho^2}.
\end{eqnarray*}
To complete the proof, it remains to take the supremum over $p_0,p_1 \in \mathfrak{D}$. \epf
\begin{cor}
\label{sor.eqat}
Let us consider $\eta \in (0,\pi/2)$ and the equation $D[t, \bldr(s)] = \gamma^\prime$, where $\bldr(s) :=\bldr+s N, N:=N(t,\bldr)$ and $t, \bldr$ are given and such that
$\gamma_-(t) < \gamma = D[t,\bldr] < \gamma_+(t)$. Whenever
$$
|\gamma - \gamma^\prime| \leq q:=\frac{\Delta_\nabla}{2}\min\big\{ \frac{\gamma - \gamma_-(t)}{b_\rho}; \frac{\gamma_+(t) - \gamma}{b_\rho}; \frac{\eta}{3 \sqrt{b_\varkappa^2+b_\tau^2}} \big\},
$$
this equation has a solution $s$ such that
$|s| \leq \frac{1}{\cos \eta \Delta_\nabla}|\gamma^\prime- \gamma|$, where $\Delta_\nabla$ is taken from \eqref{rho.lower}.
\end{cor}
\pf
We first note that due to \eqref{estim}, \eqref{rho.lower},
$\gamma_-(t) \leq D[t,\bldr(s)] \leq \gamma_+(t)$ and $|s| \leq \frac{\eta}{ \sqrt{b_\varkappa^2+b_\tau^2}}$ entail that
$$
\frac{d}{ds} D[t, \bldr(s)] = \spr{\nabla D[t,\bldr(s)]}{N} = \rho[t,\bldr(s)]
\spr{N[t,\bldr(s)]}{N} \in \left[\cos \eta \Delta_\nabla ; b_\rho \right] .
$$
So the premise of this entailment is true whenever
$$
|s| \leq \min\left\{\frac{ \eta}{ \sqrt{b_\varkappa^2+b_\tau^2}}; \frac{\gamma-\gamma_-(t)}{b_\rho}; \frac{\gamma_+(t) - \gamma}{b_\rho}\right\} = \frac{2}{\Delta_\nabla} q,
$$
and as $s$ runs over $\big[ -\frac{2}{b_\rho} q, \frac{2}{b_\rho} q\big]$, the value $D[t,\bldr(s)]$ covers a set that includes the interval $[\gamma-q, \gamma+q]$. This implies existence of solution, whereas the estimate of $|s|$ follows from the above estimate of $\frac{d}{ds} D[t, \bldr(s)]$. \epf
\par
The distance between the sets $A$ and $B$ is denoted by $\dist{A}{B} := \inf_{\bldr_a \in A, \bldr_b \in B} \|r_a-r_b\|$.
\begin{lemma}
\label{lem.isodist}
Whenever $\gamma_-(t) \leq \gamma_1 , \gamma_2 \leq \gamma_+(t)$, the following inequalities hold:
\begin{equation}
\label{isol.estm}
b_\rho^{-1}|\gamma_2-\gamma_1| \leq \dist{I(t,\gamma_1)}{I(t,\gamma_2)} \leq \Delta_\nabla^{-1} |\gamma_2-\gamma_1| .
\end{equation}
\end{lemma}
\pf
Since $t$ does not affect the claim, we will notationally ignore it in the proof.
Let $\gamma_1 \leq \gamma_2$ for the definiteness.
Since the isolines $I[\gamma_i]$ are compact due to ii) in Remark~\ref{rem.th}, there exist points $\bldr_i \in I[\gamma_i], i=1,2$ such that $\dist{I(\gamma_1)}{I(\gamma_2)} = \|\bldr_2-\bldr_1\|$. Hence for $\bldr(s):= \bldr_1 + s(\bldr_2-\bldr_1)$, we have
\begin{gather*}
\frac{d}{ds} D[\bldr(s)] = \spr{\nabla D[\bldr(s)]}{\bldr_2-\bldr_1} \leq \|\nabla D[\bldr(s)]\| \|\bldr_2-\bldr_1\| \overset{\text{\eqref{speed1}, \eqref{estim}}}{\leq} b_\rho \|\bldr_2-\bldr_1\|,
\\
\gamma_2-\gamma_1 = D[\bldr(1)] - D[\bldr(0)] =\int_0^1 \frac{d}{ds} D[\bldr(s)] \; ds \leq  b_\rho \dist{I(\gamma_1)}{I(\gamma_2)},
\end{gather*}
which implies the first inequality in \eqref{isol.estm}. For any $\gamma \in [\gamma_1, \gamma_+)$, we have by Corollary~\ref{sor.eqat},
\begin{equation}
\label{trut}
\lim_{\delta \to 0} \delta^{-1}\dist{I(\gamma+\delta)}{I(\gamma)} \leq  \Delta_\nabla^{-1}.
\end{equation}
Hence
$$
|\dist{I(\gamma_1)}{I(\gamma+\delta)} - \dist{I(\gamma_1)}{I(\gamma)}| \leq \dist{I(\gamma+\delta)}{I(\gamma)} \to 0 \quad \text{as}\quad \delta \to 0,
$$
 i.e., the map $\gamma \in [\gamma_1, \gamma_+) \mapsto \dist{I(\gamma_1)}{I(\gamma)}$ is continuous. So for any $\ve>0$, the set
$$
S:= \left\{ \gamma \in [\gamma_1, \gamma_+]: \dist{I(\gamma_1)}{I(\gamma)} \leq  (\Delta_\nabla^{-1}+\ve) (\gamma - \gamma_1)\right\}
$$
is closed and by \eqref{trut} (where $\gamma=\gamma_1$), contains an interval of the form $[\gamma_1,\gamma_\ast]$, among which there is a maximal one. For it, $\gamma_\ast = \gamma_+$. Indeed if $\gamma_\ast < \gamma_+$, \eqref{trut} (where $\gamma=\gamma_\ast$) implies that there exist $\delta_0>0$ such that $\gamma_\ast+\delta_0 \leq \gamma_+$, and for $\delta \in (0,\delta_0]$,
\begin{eqnarray*}
\frac{\dist{I(\gamma_\ast+\delta)}{I(\gamma_\ast)}}{\delta} \leq \Delta_\nabla^{-1} + \ve \Rightarrow \dist{I(\gamma_\ast+\delta)}{I(\gamma_1)} \leq \dist{I(\gamma_\ast+\delta)}{I(\gamma_\ast)} + \dist{I(\gamma_\ast)}{I(\gamma_1)}
\\
\leq (\Delta_\nabla^{-1}+\ve) \delta + (\Delta_\nabla^{-1}+\ve) (\gamma_\ast - \gamma_1) \leq (\Delta_\nabla^{-1}+\ve)(\gamma_\ast +\delta - \gamma_1) \Rightarrow
\gamma_\ast +\delta \in S,
\end{eqnarray*}
in violation of the definition of $t_\ast$ as the end of the maximal interval. This contradiction  proves that the second inequality from \eqref{isol.estm} holds for any $\gamma_2 \in [\gamma_1,\gamma_+]$, with $\Delta_\nabla^{-1}$ replaced by $\Delta_\nabla^{-1} + \ve$. The proof is completed by letting $\ve \to 0$. \epf
\begin{lemma}
\label{lem.isomov}
Suppose that $\gamma(t) \in \left[ \gamma_- (t), \gamma_+(t)\right] $ is a smooth function of time $t$ and $\bldr_\ast \in \br^2$. Then the function $t \mapsto \L(t):= \dist{I[t,\gamma(t)]}{\bldr_\ast}$ is locally Lipschitz continuous and for almost all $t$,
\begin{equation}
\label{estim.lgaga}
\left|\dot{\L}(t)\right| \leq \sup_{\bldr \in I[t,\gamma(t)]}|\lambda(t,\bldr)| + |\dot{\gamma}(t)| \sup_{\bldr \in I[t,\gamma(t)]}\rho^{-1}(t,\bldr) .
\end{equation}
\end{lemma}
\pf
Given $t = \ov{t}$, Corollary~\ref{sor.eqat} entails that $\dist{I[t,\gamma(t)]}{I[\ov{t}, \gamma(\ov{t})]} \to 0$ as $t \to \ov{t}$, and the set $\{\bldr^\prime = \bldr+s N(\ov{t}, \bldr): \bldr \in I[\ov{t}, \gamma(\ov{t})], |s| < \ve \}$ is a neighborhood of $I[\ov{t}, \gamma(\ov{t})]$ for small enough $\ve >0$. It follows that $I[t, \gamma(t)]$ lies in this neighborhood whenever $t \approx \ov{t}$, and so any point $\bldr^\prime \in I[t, \gamma(t)]$ has the form $\bldr^\prime = \bldr+s N(\ov{t}, \bldr), \bldr \in I[\ov{t}, \gamma(\ov{t})]$. On the other hand for small enough $\ve$, the implicit function theorem guarantees that $s=s(\bldr,t)$ here and moreover, the function $s(\cdot,\cdot)$ is smooth. Since
\begin{equation}
\label{lgaga}
\L(t) = \min_{\bldr \in I[\ov{t}, \gamma(\ov{t})]} \|\bldr + s(\bldr,t)N[\ov{t},\bldr ] - \bldr_\ast\|
\end{equation}
for $t \approx \ov{t}$ and the function following the $\min$ is continuous in $\bldr$ and locally Lipschitz continuous in $t$ uniformly over $\bldr \in I[\ov{t}, \gamma(\ov{t})]$, the minimum $\L(t)$ is also locally Lipschitz continuous \cite[Sect.~2.8]{Clark83} and so $\dot{\L}(t)$ exists for almost all $t$. Moreover, whenever $\dot{\L}(\ov{t})$ exists and $\ov{t} \in S:= \{\tau: D[\tau,\bldr_\ast] \neq I[\tau,\gamma(\tau)] \}$,
$$
\dot{\L}(t) = \min_{\bldr \in M} \frac{\partial}{\partial t} \|\bldr + s(\bldr,t)N[\ov{t},\bldr ] - \bldr_\ast\|\Big|_{t=\ov{t}} = \min_{\bldr \in M}
\frac{\partial s}{\partial t}(\ov{t},\bldr) \spr{N[\ov{t},\bldr]}{\frac{\bldr-\bldr_\ast}{\|\bldr-\bldr_\ast\|}},
$$
where $M$ is the set of $\bldr$'s furnishing the minimum in \eqref{lgaga}. For them
$$
\frac{\bldr-\bldr_\ast}{\|\bldr-\bldr_\ast\|} = \eta N[\ov{t},\bldr], \eta = \sgn (\gamma(\ov{t}) - D[\ov{t},\bldr_\ast]),
$$
whereas
$$
D[t, \bldr + s(\bldr,t)N(\ov{t},\bldr)] \equiv \gamma(t) \Mapsto D^\prime_t [\ov{t},\bldr] + \frac{\partial s}{\partial t}(\ov{t},\bldr) \spr{\nabla D(\ov{t},\bldr)}{N(\ov{t},\bldr)} = \dot{\gamma}(\ov{t}) \overset{\text{\eqref{speed1}}}{\Mapsto} \frac{\partial s}{\partial t}(\ov{t},\bldr) = \lambda(\ov{t},\bldr) + \frac{\dot{\gamma(\ov{t})}}{\rho(\ov{t},\bldr)},
$$
which implies \eqref{estim.lgaga} for $t = \ov{t}$ and almost all $\ov{t} \in S$. Since $\L(t) \equiv 0 \; \forall t \not\in S$, the derivative $\dot{\L}(t)$ is zero for almost all $t \not\in S$, and so  \eqref{estim.lgaga} does hold as well. \epf
\par
Now we show that the robot arrives at a state where $\dot{d}(t_\ast) \geq \nu - \delta, \sigma(t_\ast) \leq 0$ if it starts deep enough in $Z_{\text{reg}}$.
\begin{lemma}
\label{lem.start}
 Let the robot start at $t=t_0$ with $\dot{d} \leq \nu - \delta$
 and
 $$
 \gamma_-(t_0) + \Delta_\gamma/2 \leq D[t_0,\bldr(t_0)] \leq \gamma_+(t_0) -  \Delta_\gamma/2.
  $$
  Then there exists a time $t_\ast\geq t_0 $ such that
  $$
  \dot{d}(t_\ast) = \nu - \delta, \quad \ddot{d}(t_\ast) \geq 0, \quad \sigma(t_\ast) \leq 0,
  $$
  and
  $$
  [t, \bldr(t)] \in Z_{\text{reg}}, \quad  \|\bldr(t) - \bldr(t_0)\| \leq \frac{\ov{v}}{\mu \delta}(2k+1) \quad \forall t \in [t_0,t_\ast].
  $$
\end{lemma}
\pf
If $\dot{d}(t_0) = \nu - \delta, \ddot{d}(t_0) \geq 0$, the claim is clear due to \eqref{impl}.
Otherwise either $\dot{d}(t_0) < \nu - \delta$ or $\dot{d}(t_0) = \nu - \delta, \ddot{d}(t_0) < 0$; in both cases, the open set
$$
S:=\{t>t_0 :\dot{d}(t) < \nu - \delta, \gamma_-(t) < D[t,\bldr(t)] < \gamma_+(t)\}
$$
contains all $t>t_0$ that are close enough to $t_0$. Let $(t_0,t_\ast)$ be the leftmost connected component
of this set.
\par
1) Suppose that $t_\ast > t_0+ \frac{2 \pi(k+1)}{\mu \delta}$. By Lemma~\ref{lem.diberg}, there exists $T<\frac{2 \pi(k+1)}{\mu \delta}, T>0$ such that for $t\in [t_0,t_0+T]$, the robot remains in the disc $B$ of the radius $\frac{\ov{v}}{\mu \delta}(2k+1)$ centered at $\bldr(t_0)$, whereas the vector $\vec{e}(\theta)$ rotates clockwise through the angle $2 \pi (k+1)$.
By \eqref{estim}, \eqref{rho.lower}, and \eqref{estim.lgaga}, the distances from $\bldr(t_0)$ to $I[t,\gamma_i(t)]$ with $i=\pm$ are reduced by no more than
$$
(b_\lambda+\ov{\gamma}\Delta_\nabla^{-1})T < (b_\lambda+\ov{\gamma}\Delta_\nabla^{-1}) \frac{2 \pi(k+1)}{\mu \delta}.
 $$
 By \eqref{isol.estm}, these distances are initially no less than $b_\rho^{-1}\Delta_\gamma/2$.
 By the second inequality from \eqref{mu2.choice} (with $\nu:=\delta$), this means that the isolines $I[t,\gamma_-(t)]$ and $I[t,\gamma_+(t)]$ do not reach the disc $B$ and so $[t_0,t_0+T] \times B \subset Z_{\text{reg}}$. Then \eqref{speed1} and Lemma~\ref{lemrroott} ensure that as $t$ runs from $t_0$ to $t_0+T$, the gradient $\nabla D[t,\bldr(t)]$ turns through an angle that does not exceed
 $$
 b_\omega^\nabla T +  \sqrt{b_\varkappa^2+b_\tau^2} \text{\bf diam}\, B < b_\omega \frac{2 \pi(k+1)}{\mu \delta} + 2 \sqrt{b_\varkappa^2+b_\tau^2} \frac{\ov{v}}{\mu \delta}(2k+1) < 2\pi k,
 $$
 where the last inequality follows from the first relation in \eqref{mu2.choice} (with $\nu:=\delta$). Since
meanwhile $\vec{e}(\theta)$ turns through an angle of $2 \pi(k+1)$,
the gradient $\nabla D[\tau,\bldr(\tau)]$ and $\vec{e}[\theta(\tau)]$ are identically directed at some time instant $\tau \in [t_0,t_0+T]$. Then
$$
\dot{d}(\tau) = D^\prime_t[\tau, \bldr(\tau)]+ \ov{v} \spr{\nabla D[\tau, \bldr(\tau)]}{\vec{e}[\theta(\tau)]} = \rho[\tau, \bldr(\tau)] \big\{ \ov{v} - \lambda[\tau, \bldr(\tau)] \big\}  \overset{\text{\eqref{v.astchoice}}}{>} \nu,
$$
in violation of $\dot{d}(t) \leq \nu-\delta \; \forall t \in [t_0,t_\ast)$. Hence $t_\ast \leq t_0+ \frac{2 \pi(k+1)}{\mu \delta}$.
So either 2) $D[t_\ast,\bldr(t_\ast)] = \gamma_-(t_\ast)$, or 3) $D[t_\ast,\bldr(t_\ast)] = \gamma_+(t_\ast)$, or 4) $\gamma_-(t_\ast) < D[t_\ast,\bldr(t_\ast)] < \gamma_+(t_\ast)$ and $\dot{d}(t_\ast)=\nu - \delta$.
\par
2,3) By retracing the above arguments, we see that in case 2), the isoline cannot reach the disk $B$ and thus the current location of the robot in time $t_\ast-t_0 \leq \frac{2 \pi(k+1)}{\mu \delta}$. Thus case 2) does not occur. Case 3) is similarly rejected.
\par
4) In this case, $\dot{d}(t) \leq \nu - \delta \; \forall t \in [t_0,t_\ast] \wedge \dot{d}(t_\ast) = \nu - \delta \Mapsto \ddot{d}(t_\ast) \geq 0 \Mapsto \sigma(t_\ast) \leq 0$, where the last implication is due to \eqref{impl}. The last claim of the lemma follows from Lemma~\ref{lem.diberg}. \epf
\par
By summarizing the foregoing, now we show that the robot does reach $V_\star(t)$ with a certain ``excess''.
\begin{lemma}
\label{lem.reach}
In a finite time, the robot reaches the set $V_{\star\star}(t):=\{\bldr : D(t,\bldr) > \gamma_{\star\star}(t)\}$, where $\gamma_{\star\star}(t):= \gamma_\star(t)+\Delta_\gamma/2$.
\end{lemma}
\pf Suppose to the contrary that the claim is incorrect
$
\bldr(t) \not\in V_{\star\star}(t) \Leftrightarrow D[t,\bldr(t)] \leq \gamma_{\star\star} (t)\; \forall t \geq 0
$,
where $\gamma_{\star\star}(t) \leq \gamma_+(t) - \Delta_\gamma/2$ by \eqref{rho.lower}. We consider separately three cases.
\par
{\bf (a)} $\nu - \delta \leq \dot{d}(0), \sigma(0) \leq 0$.
This case is impossible due to Corollary~\ref{cor.mcase}.
\par
{\bf (b)} $\nu - \delta \geq \dot{d}(0)$. By Lemma~\ref{lem.start}, along with \eqref{rho.lower}, this case entails (a) at some time instant and so is impossible.
\par
{\bf (c)} $\dot{d}(0) > \nu - \delta$. Let $[0,t_\ast)$ be the leftmost connected component of the set
$$
\{t \geq 0 : \dot{d}(t) > \nu - \delta, \gamma_- < D[t,\bldr(t)] < \gamma_+(t)\}.
 $$
 Since $\nu - \delta >0$ and by \eqref{estim}, $\infty>\gamma_+^0 \geq \gamma_+(t) - \Delta-\gamma/2 \geq \gamma_{\star\star}(t) \geq D[t,\bldr(t)]$, we conclude that $t_\ast < \infty$. Since $\dot{d}(t) > \nu - \delta > \ov{\gamma}$ for $t \in (0,t_\ast)$ by \eqref{delta.choice} and $|\dot{\gamma}_-(t)| \leq \ov{\gamma}$, we see that $\dot{d}(t_\ast) = \nu-\delta$. The proof is completed by retracing the arguments from b).
\par
The contradictions obtained complete the proof.
\epf
\par
The last lemma in fact completes the proofs of Theorems~\ref{th.m0}, \ref{th.m1} and Remark~\ref{rem.th}.
\begin{lemma}
\label{lem.keep}
The robot cannot leave the desired vicinity $V_\star(t)$ of the maximizer when starting at $t=\tau$ in $V_{\star\star}(\tau)$.
\end{lemma}
\pf
Suppose to the contrary that the robot leaves $V_\star(t)$, starting at $t=\tau$ with $\bldr(\tau) \in V_{\star\star}(\tau) \Leftrightarrow D[\tau,\bldr(\tau)] > \gamma_{\star\star}(\tau)$. By the continuity argument, it intersects first $I[t,\gamma_{\star\star}(t)]$ and then $I[t,\gamma_\star(t)]$
Let $t_1>0$ be the earliest time $t$ such that $d(t) = \gamma_\star(t)$, and let $t_0$ be the latest time $t \in (0,t_1)$ such that $d(t) = \gamma_{\star\star}(t)$.
Then
$$
d(t) < \gamma_{\star\star}(t) \; \forall t \in (t_0,t_1) \Rightarrow \dot{d}(t_0) \leq \dot{\gamma}_{\star\star}(t_0) = \dot{\gamma}_{\star}(t_0) \leq \nu - \delta
$$
by \eqref{estim} and \eqref{delta.choice}.
By Lemma~\ref{lem.start}, there is $t_\ast\geq t_0 $ such that
$\dot{d}(t_\ast) = \nu - \delta, \sigma(t_\ast) \leq 0$, and $\|\bldr(t) - \bldr(t_0)\| \leq \frac{\ov{v}}{\mu \delta}(2k+1) \; \forall t \in [t_0,t_\ast]$. This inequality and the arguments from 1) in the proof of Lemma~\ref{lem.start} show that
the robot does not reach $I[t,\gamma_\star(t)]$ for $t \in [t_0,t_\ast]$ and so $\gamma_\star(t_\ast) < d(t) \leq \gamma_{\star\star}(t_\ast)$.
Then for $t \geq t_\ast$, the robot reaches $I[t,\gamma_+(t)]$ with constantly maintaining $\dot{d}\geq \nu - \delta$ true, where $\nu - \delta > \dot{\gamma}_\star(t)$. So $d(t), t \geq t_0$ cannot reach the value $\gamma_\star(t)$ earlier than it reaches $\gamma_{\star\star}(t)$ for the second time, in violation of the definitions of $t_0,t_1$. The contradiction obtained proves the lemma. \epf
\par
{\bf \uppercase{proof of theorem}~\ref{th.m1}:} This theorem is immediate from Lemmas~\ref{lem.reach} and \ref{lem.keep}. \epf
\par
{\bf \uppercase{proof of theorem}~\ref{th.m0}:} This theorem is immediate from Theorem~\ref{th.m1}. \epf
\par
{\bf \uppercase{proof of remark}~\ref{rem.th}:} This remark holds since in this appendix, Assumption~\ref{ass.0a} was taken in the generalized form of i)---iii) from the remark. \epf

\bibliographystyle{plain}
  \bibliography{Hamidref}
 \end{document}